\documentclass[a4paper,final]{article}

\usepackage{amsthm}
\newtheorem{teo}{Theorem}[section]
\newtheorem{lemma}[teo]{Lemma}
\newtheorem{prop}[teo]{Proposition}
\newtheorem{cor}[teo]{Corollary}

\theoremstyle{definition}

\newtheorem{rem}[teo]{Remark}

\newcommand{\dimo}[1]{\vspace{1pt}\noindent{\it Proof of} {\hspace{2pt}}\ref{#1}.}
\newcommand{\finedimo}[1]{{\hfill\hbox{\enspace\fbox{\ref{#1}}}}\vspace{5pt}}

\usepackage{amssymb}
\usepackage{amsmath}

\usepackage[]{graphicx}

\newcommand{\MPa}{\mbox{${\rm MPa}$}}
\newcommand{\Va}{\mbox{${\rm Va}$}}
\newcommand{\La}{\mbox{${\rm La}$}}
\newcommand{\Ba}{\mbox{${\rm Ba}$}}
\newcommand{\Ca}{\mbox{${\rm Ca}$}}
\newcommand{\calT}{{\cal T}}
\newcommand{\calP}{{\cal P}}
\newcommand{\calQ}{{\cal Q}}
\newcommand{\calX}{{\cal X}}
\newcommand{\calB}{{\cal B}}
\newcommand{\calN}{{\cal N}}
\newcommand{\calL}{{\cal L}}
\newcommand{\st}[1]{{\rm star}(#1)}
\newcommand{\clst}[1]{{\rm clst}(#1)}
\newcommand{\inter}[1]{{\rm Int}(#1)}

\newcommand{\MPag}{\mbox{${\rm MPd}$}}
\newcommand{\Lag}{\mbox{${\rm Ld}$}}
\newcommand{\Mag}{\mbox{${\rm Md}$}}
\newcommand{\Nag}{\mbox{${\rm Nd}$}}
\newcommand{\swell}[1]{{\rm sw}(#1)}

\title{A calculus for ideal triangulations of three-manifolds with
  embedded arcs}

\author{Gennaro Amendola}

\begin{document}

\maketitle

{\small\noindent{\sc Abstract}:
Refining the notion of an ideal triangulation of a compact
three-manifold, we provide in this paper a combinatorial presentation
of the set of pairs $(M,\alpha)$, where $M$ is a three-manifold and
$\alpha$ is a collection of properly embedded arcs.
We also show that certain well-understood combinatorial moves are
sufficient to relate to each other any two refined triangulations
representing the same $(M,\alpha)$.
Our proof does not assume the Matveev-Pergallini calculus for ideal
triangulations, and actually easily implies this calculus.}

\vspace{.5cm}


\vspace{.5cm}


\section*{Introduction}

A {\em combinatorial presentation} of a class of topological objects
(viewed up to the appropriate equivalence relation) is a set of finite
combinatorial objects, such that each combinatorial object defines
(say ``presents'') a unique topological object, and each topological
object is presented by at least one combinatorial object.
A {\em calculus} for a combinatorial presentation is a finite set of
moves on the combinatorial objects, such that two combinatorial
objects present the same topological object if and only if they are
related to each other by a finite sequence of moves in the given set.

Combinatorial presentations are fundamental tools for studying
3-manifolds and links, and for constructing invariants.
They translate a topological problem into a combinatorial and, maybe,
a simpler one.
For instance, an invariant on the class of topological objects can be
defined on the combinatorial objects, checking that it is preserved by
the moves.

For 3-manifolds, there are several different types of presentations,
{\em e.g.}~(ideal) triangulations, Heegaard diagrams, surgery (on
links), and spines.
In the present work we concentrate on the pairs $(M,\alpha)$, where
$M$ is a compact connected $3$-manifold with non-empty boundary and
$\alpha = \{\alpha^{(1)}, \ldots, \alpha^{(n)}\}$ is a (possibly
empty) collection of disjoint arcs properly embedded in $M$ (viewed up
to simultaneous isotopy).
We provide a presentation of such pairs and we describe the
corresponding calculus.
The objects of the presentation are the {\em marked ideal
triangulations} of the pair
$(M,\alpha)$, that is the ideal triangulations of $M$ that contain as edges
all the arcs in $\alpha$, and the moves of the calculus are the moves
on ideal triangulations ({\em i.e.}~Matveev-Piergallini moves) which
do not kill edges belonging to $\alpha$ (such moves will be called
{\em admissible}).

The calculus for marked ideal triangulations is not new: in fact it has
been used by Baseilhac and Benedetti (see~\cite{BB1,BB2,BB3}) in the prove
that the so-called {\it quantum hyperbolic invariants} (QHI) for
links in 3-manifolds equipped with flat $PSL(2,\mathbb{C})$-bundles
are well defined.
They derived this calculus from the Matveev-Piergallini
one~\cite{Matveev:calculus,Piergallini}, as refined by Turaev and Viro
in~\cite{Turaev-Viro}.
They have also used the generalization to the setting of marked ideal
triangulations of a result of Makovetskii~\cite{makov}.
We will give a new proof of the calculus (for marked ideal
triangulations), which is instead self-contained, see
Section~\ref{calculus:sec}.
Actually, our proof specializes to a new proof of the
Matveev-Piergallini calculus.
Although our proof is quite long, it is conceptually very simple: in
fact it uses only easy results on triangulations and easy topological
arguments.
For the sake of completeness, we will also describe a sketch of the
derivation of the calculus for marked triangulations from the
Matveev-Piergallini one, see Subsection 2.5.

The generalized Makovetskii result states that, if two marked ideal
triangulations of a pair $(M,\alpha)$ are given, then they are
dominated, as far as some {\em positive} admissible moves are
concerned, by another marked ideal triangulation.
An admissible move is positive if it increases the number of tetrahedra.
In Section~\ref{dominating:sec}, we provide the details of the proof
of this refinement, and, in Subsection~\ref{hamiltonian:subsec}, we
describe the relationship between marked ideal triangulations and
links in 3-manifolds.

The initial motivation of the present paper was the remark, due to
Frigerio and Petronio~\cite{Frig-Petr}, that marked ideal
triangulations naturally arise in the study of {\em complete
  finite-volume orientable hyperbolic $3$-manifolds with geodesic
  boundary.}
In Subsection~\ref{part_trunc_tria:subsec} we will describe how this
relationship arises.

\section{Definitions}

From now on, unless explicitly stated, $M$ will be a compact connected
$3$-manifold with non-empty boundary, and $\alpha = \{\alpha^{(1)},
\ldots, \alpha^{(n)}\}$ will be a (possibly empty) collection of
disjoint arcs properly embedded in $M$, viewed up to simultaneous
isotopy.

\subsection{Standard spines and moves}\label{spines_moves:subsection}

In this subsection we recall the definition of spine and we describe
some moves.

\paragraph{Standard spines}
A {\em quasi-standard} polyhedron $P$ is a finite, connected, and
purely $2$-dimensional polyhedron with singularities of stable nature
({\em i.e.}~triple lines and points where 6 non-singular components
meet).
Such a polyhedron is called {\em standard} if it is cellularized by
singularity (depending on dimension, we call the components {\em
vertices}, {\em edges}, and {\em regions}).
A quasi-standard sub-polyhedron $P$ of $M$ contained in $\inter{M}$ is
called a {\em spine} of $M$ if the manifold $M$ collapses to it (or,
equivalently, $M \setminus P \cong \partial M \times [0,1)$).
Each spine of $M$ is always viewed up to isotopy.
For the sake of completeness, let us recall that, if $M$ is closed,
the boundary is created by puncturing $M$ ({\em i.e.}~by considering
$M$ minus a ball).

It is by now well-known, after the work of Casler~\cite{Casler}, that
a standard spine determines $M$ uniquely up to homeomorphism and that
every $M$ has standard spines.
In the sequel we will omit the word ``standard'', writing only
``spine''; nevertheless, if standardness will not be obvious, we will
use the word ``standard''.
Moreover, in the figure of a piece of spine the singular set is drawn
thick.

\begin{figure}
  \centerline{\includegraphics{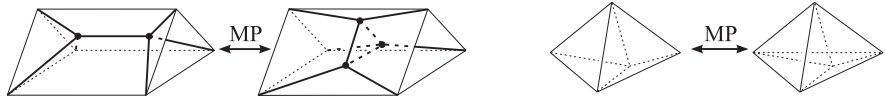}}
  \caption{The MP-move on a spine (left) and on the dual ideal
  triangulation (right).}
  \label{MP_move:fig}
\end{figure}

\paragraph{MP-move}
Any two (standard) spines of $M$ can be transformed into each other by
certain well-understood moves.
Let us start from the move shown in Fig.~\ref{MP_move:fig}-left, which
is called MP-{\em move}.
Such a move will be called {\em positive} if it increases (by one) the
number of vertices, and {\em negative} otherwise.
Note that, if we apply an MP-move to a spine of $M$, the result will
be another spine of $M$.
It is already known (but it will also follow from our
Corollary~\ref{gener_MP_senza_V:cor}), after the work of
Matveev~\cite{Matveev:calculus} and Piergallini~\cite{Piergallini},
that any two standard spines of the same $M$ with at least two
vertices can be transformed into each other by MP-moves (see
Theorem~\ref{MP_calculus:teo}).

\begin{figure}
  \centerline{\includegraphics{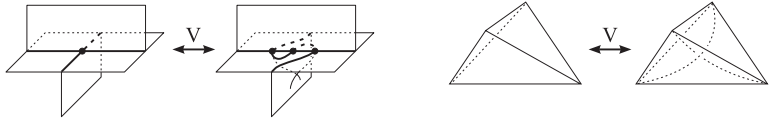}}
  \caption{The V-move on a spine (left) and on the dual ideal
  triangulation (right).}
  \label{V_move:fig}
\end{figure}

\paragraph{V-move}
If one of the two spines of $M$ (we want to transform into each other)
has just one vertex, another move is required.
The move shown in Fig.~\ref{V_move:fig}-left is called V-{\em move}.
Note that if we apply such a move to a spine of $M$, the result will
be another spine of $M$.
As above, we have {\em positive} and {\em negative} V-moves.
Note that $3$ different positive V-moves can be applied at each vertex.

\begin{figure}
  \centerline{\includegraphics{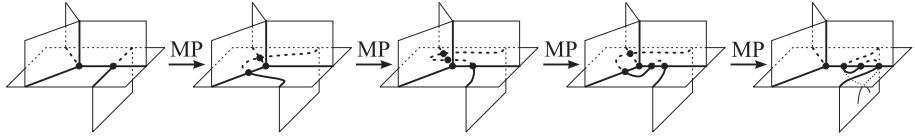}}
  \caption{If there is another vertex, each positive V-move is a
  composition of MP-moves.}
  \label{V_comp_MP:fig}
\end{figure}

If a positive V-move is applied to a spine with at least two vertices,
the V-move is a composition of MP-moves.
In Fig.~\ref{V_comp_MP:fig} we show the three positive and the one
negative MP-moves giving the V-move.

\begin{figure}
  \centerline{\includegraphics{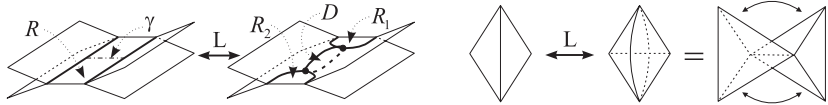}}
  \caption{The L-move on a spine (left) and on the dual ideal
    triangulation (right).}
  \label{L_move:fig}
\end{figure}

\paragraph{L-move}
A generalization of the V-move is the L-{\em move}, see
Fig.~\ref{L_move:fig}-left.
As above, we have {\em positive} and {\em negative} L-moves.
As opposed to the V-move, this move is non-local, so it must be
described with some care.
A positive L-move, which increases by two the number of vertices, is
determined by an arc $\gamma$ properly embedded in a region $R$ of $P$.
The move acts on $P$ as in Fig.~\ref{L_move:fig}-left, but, to define
its effect non-ambiguously, we must specify which pairs of regions,
out of the four regions incident to $R$ at the endpoints of $\gamma$,
will become adjacent to each other after the move.
This is achieved by noting that $R$ is a disc, so its regular
neighborhood in $M$ is a product, and we can choose for $R$ a
transverse orientation.
Using it, at each endpoint of $\gamma$ we can tell from each other the
two regions incident to $R$ as being an upper and a lower one, and we
can stipulate that the two upper regions will become incident after
the move (and similarly for the lower ones).
Obviously, a positive L-move leads to a (standard) spine $P'$ of $M$.

For the negative case the situation is more complicated.
A negative L-move can lead to a non-standard spine.
If $R_1$ and $R_2$ are contained in the same region, after the
negative L-move, the ``region'' $R$ would not be a disc.
To avoid this loss of standardness, we will call negative L-moves only
those preserving standardness.
So a negative L-move can be applied only if the regions $R_1$ and
$R_2$ are different.
With this convention, if we apply an L-move to a spine of $M$, the
result will be another spine of $M$.

\begin{figure}
  \centerline{\includegraphics{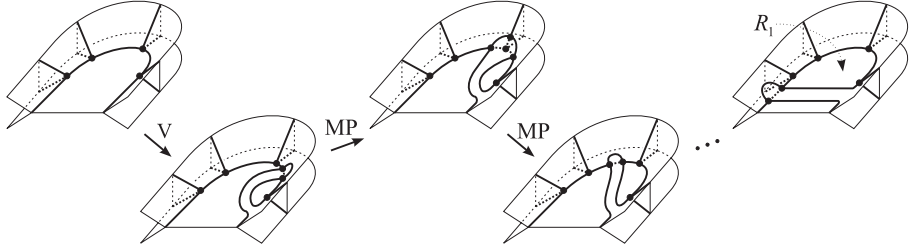}}
  \caption{Each positive L-move is a composition of V-~and MP-moves
  (case where $R_1$ has more than one vertex).}
  \label{L_comp_V_MP_1:fig}
\end{figure}

\begin{figure}
  \centerline{\includegraphics{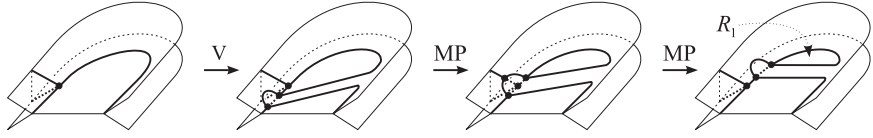}}
  \caption{Each positive L-move is a composition of V-~and MP-moves
  (case where $R_1$ has only one vertex).}
  \label{L_comp_V_MP_2:fig}
\end{figure}

Each positive L-move is a composition of V-~and MP-moves.
In Fig.~\ref{L_comp_V_MP_1:fig} we show the one positive V-move and
the pairs of (one positive and one negative) MP-moves giving the
L-move.
Obviously, to apply such moves, $R_1$ must have at least two vertices.
If $R_1$ has only one vertex, then $R_2$ has at least two vertices
(because $P'$ is standard); so we can take the symmetric picture.
For future reference, we note that, if $R_1$ has only one vertex, we
can obtain the L-move also as a composition of only one V-~and one
pair of MP-moves, as shown in Fig.~\ref{L_comp_V_MP_2:fig}.

\begin{figure}
  \centerline{\includegraphics{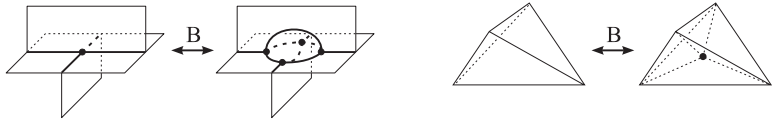}}
  \caption{The B-move on a spine (left) and on the dual ideal
  triangulation (right).}
  \label{B_move:fig}
\end{figure}

\paragraph{B-move}
Now we describe the B-{\em move} (shown in
Fig.~\ref{B_move:fig}-left).
As above, we have {\em positive} and {\em negative} B-moves.
This move is quite different from the previous ones, because if we
apply a positive B-move to a spine $P$ of $M$, the result will be a
spine $P_B$ of $M \setminus B^3$ (where $B^3$ is a 3-ball with closure
embedded in $M$).
So it is obvious that a B-move cannot be a composition of V-~and
MP-moves.
By definition of spine, we have that $M \setminus P_B$ is the disjoint
union of $\partial M \times [0,1)$ and $B^3 \cup (\partial B^3 \times
[0,1))$.
The ball $\calB = B^3 \cup (\partial B^3 \times [0,1))$ will be called
{\em proper ball}.

\begin{figure}
  \centerline{\includegraphics{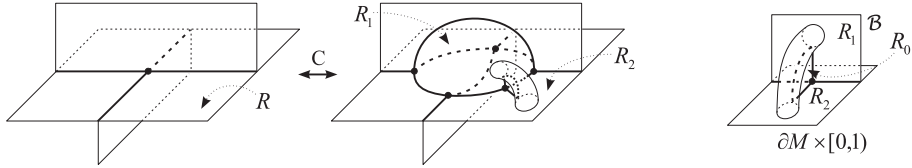}}
  \caption{The C-move on a spine (left) and the corresponding arch
  (right).}
  \label{C_move:fig}
\end{figure}

\begin{figure}
  \centerline{\includegraphics{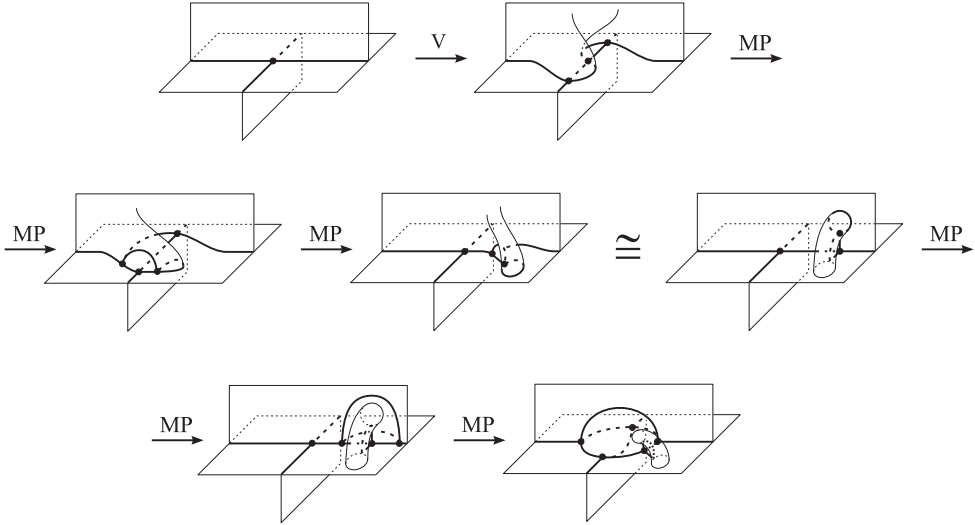}}
  \caption{Each positive C-move is a composition of V-~and MP-moves.}
  \label{C_comp_V_MP:fig}
\end{figure}

\paragraph{C-move}
In the end, we describe the C-{\em move}, see
Fig.~\ref{C_move:fig}-left.
As above, we have {\em positive} and {\em negative} C-moves.
This move is very similar to the B-move, but, if we apply a positive
C-move to a spine of $M$, we obtain another spine of the same $M$.
In fact, each positive C-move is a composition of V-~and MP-moves: the
V-move and the (four) MP-moves are shown in
Fig.~\ref{C_comp_V_MP:fig}.
Note also that $12$ different positive C-moves can be applied at each
vertex.

We will call {\em arch} the configuration shown in
Fig.~\ref{C_move:fig}-right, created by a C-move.
Let us compare the spine $P_C$, obtained from a spine $P$ via a
C-move, with the spine $P_B$, obtained from $P$ via a B-move (applied
at the same vertex).
They are different only for the presence of the arch, which joins two
different regions ($R_1$ and $R_2$) of the spine $P_B$.
Note also that, after the C-move, the proper ball $\calB$, created by
the B-move, is connected to $\partial M \times [0,1)$ by the cavity of
the arch.

\paragraph{}
In the rest of the paper we will always regard $M$ as being fixed and
we will only consider spines and moves embedded in $M$, without
explicit mention.

\subsection{Ideal triangulations}\label{id_tria:subsection}

In this subsection we recall the definition of loose triangulation and
ideal triangulation, eventually defining the marked ideal
triangulations (and spines) of a pair $(M,\alpha)$ and the moves on
them.

\paragraph{Loose and ideal triangulations}
A {\em loose triangulation} of a polyhedron $|\calP|$ is a
triangulation $\calP$ of $|\calP|$ in a weak sense, namely
self-adjacencies and multiple adjacencies are allowed.
For any manifold $M$ (as above), let us denote by $\widehat{M}$ the
space obtained from $M$ by collapsing to a point each component of
$\partial M$.
An {\em ideal triangulation} of a manifold $M$ (as above) is a
partition $\calT$ of $\inter{M}$ into open cells of dimensions 1, 2,
and 3, induced by a loose triangulation $\widehat{\calT}$ of the space
$\widehat{M}$ such that the vertices of $\widehat{\calT}$ are
precisely the points of $\widehat{M}$ corresponding to the components
of $\partial M$.
The quotient of $\partial M$ will be denoted by $\widehat{\partial
M}$.
Note that $\widehat{M} \setminus \widehat{\partial M}$ can be
identified with $\inter{M}$.
As for spines, each ideal triangulation of $M$ is always viewed up to
isotopy.

\begin{figure}
  \centerline{\includegraphics{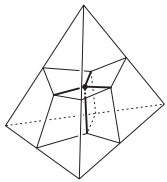}}
  \caption{Portion of spine dual to a tetrahedron of an ideal
  triangulation.}
  \label{duality:fig}
\end{figure}

\paragraph{Duality}
We show now the well-known fact that ideal triangulations exist for
each $M$.
It turns out \cite{Matveev-Fomenko,Petronio:tesi,Matveev:new:book}
that there exists a natural bijection between standard spines and
ideal triangulations of a 3-manifold.
Given an ideal triangulation $\calT$, the corresponding standard spine
$P$ is just the 2-skeleton of the dual cellularization, as illustrated
in Fig.~\ref{duality:fig}.
The inverse passage is also explicit, but it is a little more
difficult; so we omit its description.
The ideal triangulation $\calT$ and the spine $P$ are said to be {\em
dual}.
As said above, every $M$ has standard spines, so dually it has ideal
triangulations.

We show in
Figg.~\ref{MP_move:fig}-right,~\ref{V_move:fig}-right,~\ref{L_move:fig}-right,
and~\ref{B_move:fig}-right the MP-,~V-,~L-, and~B-moves, respectively,
on a spine in terms of the dual ideal triangulations (we have omitted
the dual version of the C-move because of the complexity of the
picture).
In the sequel we will intermingle the spine and the ideal
triangulation viewpoints.

\paragraph{Marked ideal triangulations (and spines)}
Recall that $\alpha$ is a collection of disjoint arcs properly
embedded in a manifold $M$.
A {\em marked ideal triangulation} of the pair $(M,\alpha)$ is a pair
$(\calT ,\beta)$, where $\calT$ is an ideal triangulation of $M$ and
$\beta = \{\beta^{(1)},\ldots ,\beta^{(n)}\}$ is a collection of edges
of $\calT$ (simultaneously) isotopic to $\alpha =
\{\alpha^{(1)},\ldots ,\alpha^{(n)}\}$.
The quotient of $\beta = \{\beta^{(1)},\ldots ,\beta^{(n)}\}$ in
$\widehat{\calT}$ will be denoted by $\widehat{\beta} =
\{\widehat{\beta}^{(1)},\ldots ,\widehat{\beta}^{(n)}\}$, and the pair
$(\widehat{\calT},\widehat{\beta})$ will be said {\em marked loose
triangulation corresponding to $(\calT ,\beta)$}.
With a little abuse of terminology, in the sequel we will say that the
edges in $\beta$ and $\widehat{\beta}$ {\em belong to $\alpha$}.

Using duality, we can give a natural definition of {\em marked spine}
of a pair $(M,\alpha)$ as a pair $(P,\tilde\beta)$, where $P$ is the
spine dual to a marked ideal triangulation $(\calT ,\beta)$ of the
pair $(M,\alpha)$ and $\tilde\beta = \{\tilde\beta^{(1)},\ldots
,\tilde\beta^{(n)}\}$ is the collection of the regions of $P$ dual to
the $\beta^{(i)}$'s.
With a little abuse of notation, we will drop the tilde, writing only
$\beta^{(i)}$ instead of $\tilde\beta^{(i)}$, and we will say that the
regions $\beta^{(i)}$ also {\em belong to $\alpha$}.

\paragraph{Existence of marked ideal triangulations}
By duality, to prove that every pair $(M,\alpha)$ has marked ideal
triangulations, it is enough to prove that it has marked ideal
spines.
So we prove that $M$ has a spine such that $\alpha$ is isotopic to the
collection of the edges dual to $n$ different regions.
Let $N(\alpha) = \sqcup_{i=1}^n N(\alpha^{(i)})$ be a regular
neighborhood of $\alpha$, let $Q$ be a spine of $M \setminus
N(\alpha)$.
Note that we have a retraction $\pi$ of $M \setminus N(\alpha)$ onto
$Q$.
For $i=1,\ldots ,n$, let $D^{(i)}$ be a 2-disc properly embedded in
$N(\alpha^{(i)})$, embedded in $\inter{M}$, and intersecting
$\alpha^{(i)}$ transversely in one point.
Now, we can suppose that, by projecting the $\partial D^{(i)}$'s to
$Q$ along $\pi$, we obtain ``half-open'' annuli $\partial D^{(i)}
\times [0,1)$.
Up to isotopy, we can also suppose that each $\pi(\partial D^{(i)})$
intersects the singularity of $P$, and that $\pi(\cup_{i=1}^n \partial
D^{(i)})$ is transversal to the singularity and to itself.
Let us define $P$ as the union of the polyhedron $Q$, the discs
$D^{(i)}$, and the annuli $\partial D^{(i)} \times [0,1)$.
Obviously, $P$ is the desired spine: in fact, $P$ is a (standard)
spine of $M$ and each $\alpha^{(i)}$ coincides with the edge dual to
the region $D^{(i)} \cup (\partial D^{(i)} \times [0,1))$ of $P$.

\paragraph{Admissible moves}
We will now discuss an extension of the MP-,~V-,~L-,~and C-moves to
the context of marked ideal triangualations.
Given a marked ideal triangulation $(\calT,\beta)$ of $(M,\alpha)$,
the idea is to consider a move from $\calT$ to $\calT'$ {\em
admissible} if there is a $\beta'$ such that $(\calT',\beta')$ is a
marked ideal triangulation of $(M,\alpha)$, and $\beta'$ coincides
with $\beta$ except ``near'' the portion of $\calT$ affected by the
move.
As it turns out, admissibility depends on $\beta$.
Moreover, $\beta'$ is sometimes not unique.

By duality, we describe the moves on spines to refer to simpler
pictures, but we invite the reader to figure out the dual ideal
triangulation pictures.
Let $(P,\beta)$ be the marked spine dual to $(\calT,\beta)$.
We describe the moves one by one.

\subparagraph{MPa-move}
A positive MP-move from $P$ to $P'$ is admissible whatever $\beta$,
and $\beta'$ consists of the same regions as $\beta$ ({\em i.e.}~the
newborn triangular region does not belong to $\beta'$); the move from
$(P,\beta)$ to $(P',\beta')$ is called {\em positive \MPa-move}.
A negative MP-move from $P$ to $P'$ is admissible if it is the inverse
of a positive \MPa-move: namely, the triangular region disappearing
during the move must not belong to $\beta$, and $\beta'$ consists of
the same regions as $\beta$; the move from $(P,\beta)$ to
$(P',\beta')$ is called {\em negative \MPa-move}.
See Fig.~\ref{MP_move:fig}.

\subparagraph{Va-move}
A positive V-move from $P$ to $P'$ is admissible whatever $\beta$, and
$\beta'$ consists of the same regions as $\beta$ ({\em i.e.}~the two
newborn little regions do not belong to $\beta'$); the move from
$(P,\beta)$ to $(P',\beta')$ is called {\em positive \Va-move}.
A negative V-move from $P$ to $P'$ is admissible if it is the inverse
of a positive \Va-move: namely, the two little regions disappearing
during the move must not belong to $\beta$, and $\beta'$ consists of
the same regions as $\beta$; the move from $(P,\beta)$ to
$(P',\beta')$ is called {\em negative \Va-move}.
See Fig.~\ref{V_move:fig}.

Now, recall that, if there are at least two vertices, a positive
V-move is a composition of MP-moves, see Fig.~\ref{V_comp_MP:fig}; the
``admissible'' version of this fact is not so obvious but it is true.
Namely, if there are at least two vertices, a positive \Va-move is a
composition of \MPa-moves.
To prove this, it is enough to note that a positive \Va-move is a
composition of MP-moves (see again Fig.~\ref{V_comp_MP:fig}), that the
negative MP-move of the sequence eliminates a region created by a
previous positive \MPa-move (so the region does not belongs to
$\beta$), and that the position of the $(\beta')^{(i)}$'s after the
\MPa-moves is the same as after the \Va-move.

\subparagraph{La-move}
For the L-moves, the situation is more complicated.
A positive L-move from $P$ to $P'$ is admissible whatever $\beta$, but
$\beta'$ is not uniquely determined.
We follow the notation of Fig.~\ref{L_move:fig}.
We have two cases depending on whether $R$ belongs to $\beta$ or not.
If $R$ does not belong to $\beta$, then $\beta'$ consists of the same
regions as $\beta$ ({\em i.e.}~$R_1$, $R_2$, and the newborn little
region $D$ do not belong to $\beta'$).
In such a case, the move from $(P,\beta)$ to $(P',\beta')$ is called
{\em positive \La-move}.
If $R$ belongs to $\beta$, the situation is a little ambiguous: $R$ is
divided in two regions, and both of them ``are isotopic to $R$'' ({\em
i.e.}~the dual edges of $R_1$ and $R_2$ are both isotopic to the dual
edge of $R$).
If we define $\beta'_1$ as $(\beta \setminus \{R\}) \cup \{R_1\}$ and
$\beta'_2$ as $(\beta \setminus \{R\}) \cup \{R_2\}$, we have two
admissible L-moves underlying the original L-move: one from
$(P,\beta)$ to $(P',\beta'_1)$ and one from $(P,\beta)$ to
$(P',\beta'_2)$.
Also both these moves are called {\em positive \La-moves}.
Note that the choice of the region, between $R_1$ and $R_2$, is
included in the move.

A negative L-move from $P$ to $P'$ is admissible if it is the inverse
of a positive \La-move.
Necessarily, the little region $D$ disappearing during the move must
not belong to $\beta$, and only one region between $R_1$ and $R_2$ can
belong to $\beta$.
Now, we have two cases: if both $R_1$ and $R_2$ do not belong to
$\beta$, then $\beta'$ consists of the same regions as $\beta$;
otherwise, if one region $R_i$ (between $R_1$ and $R_2$) belongs to
$\beta$, then $\beta'$ is equal to $(\beta \setminus \{R_i\}) \cup
\{R\}$.
In both cases, the move from $(P,\beta)$ to $(P',\beta')$ is called
{\em negative \La-move}.

Now, recall that each L-move is a composition of V-~and MP-moves (see
Figg.~\ref{L_comp_V_MP_1:fig} and~\ref{L_comp_V_MP_2:fig}).
As above, we show that each positive \La-move is a composition of
\Va-~and \MPa-moves.
If $R$ does not belong to $\beta$, the situation is analogous to that
of \Va-move, so we omit its treatment.
On the contrary, we suppose that $R$ belongs to $\beta$.
Now, one region, between $R_1$ and $R_2$, belongs to $\beta'$: we
suppose that $R_2$ belongs to $\beta'$ (the case for $R_1$ is
symmetric).
The V-~and MP-moves shown in Figg.~\ref{L_comp_V_MP_1:fig}
and~\ref{L_comp_V_MP_2:fig} (we have two cases depending on whether
$R_1$ has one vertex or more) are all admissible, and the (positive)
\Va-move leaves just $R_2$ in $\beta'$.

\subparagraph{Ca-move}
A positive C-move from $P$ to $P'$ is admissible whatever $\beta$, and
$\beta'$ consists of the same regions as $\beta$ ({\em i.e.}~the four
newborn regions do not belong to $\beta'$); the move from $(P,\beta)$
to $(P',\beta')$ is called {\em positive \Ca-move}.
See Fig.~\ref{C_move:fig}.
Note that $R_1$ is joined to $R_2$, so, if $R$ belongs to $\beta$,
then the region containing $R_1$ and $R_2$ belongs to $\beta'$.

A negative C-move from $P$ to $P'$ is admissible if it is the inverse
of a positive \Ca-move: namely, the four regions (included the disc of
the arch) disappearing during the move must not belong to $\beta$, and
$\beta'$ consists of the same regions as $\beta$.
The move from $(P,\beta)$ to $(P',\beta')$ is called {\em negative
\Ca-move}.

As above, it is easy to see that each \Ca-move is a composition of
\Va-~and \MPa-moves.

\subparagraph{Ba-move}
For the B-moves, the situation is quite different because such moves
change the homeomorphism class of the manifold.
A B-move from $P$ to $P'$ will be considered {\em admissible} both if
it is positive, or if it is negative and the four regions disappearing
do not belong to $\beta$.
In such a case, $\beta'$ consists of the same regions as $\beta$.
The move from $(P,\beta)$ to $(P',\beta')$ is called {\em \Ba-move}
({\em positive} or {\em negative}, respectively).
See Fig.~\ref{B_move:fig}.

\paragraph{}
From now on, since a marked ideal triangulation $(\calT,\beta)$ is a
pair while an ideal triangulation $\calT$ is not, then, for the sake
of shortness, we will omit the word ``marked'' (also for spines and
loose triangulations) unless the difference is not clear.

\section{The calculus}\label{calculus:sec}

The main result of this paper is the following.

\begin{teo}\label{gener_MP:teo}
Two marked ideal triangulations of a pair $(M,\alpha)$ can be obtained
from each other via a sequence of \Va-~and \MPa-moves.
\end{teo}

Recalling that, if there are at least two tetrahedra, each \Va-move is
a composition of \MPa-moves, we obtain the following corollary of
Theorem~\ref{gener_MP:teo}.

\begin{cor}\label{gener_MP_senza_V:cor}
Two marked ideal triangulations of $(M,\alpha)$ with at least two
tetrahedra can be obtained from each other via a sequence of
\MPa-moves only.
\end{cor}

As a particular case we obtain the Matveev-Piergallini theorem.

\begin{teo}[Matveev-Piergallini]\label{MP_calculus:teo}
Two spines of $M$ can be obtained from each other via a sequence of
{\rm V}-~and {\rm MP}-moves.
If moreover both spines have at least two vertices, then they can
be obtained from each other via a sequence of {\rm MP}-moves only.
\end{teo}

The idea of the proof of Theorem~\ref{gener_MP:teo} consists of the
following steps:
\begin{itemize}

\item a ``desingularization'' of the two marked ideal triangulations,
  say $(\calT_1,\beta_1)$ and $(\calT_2,\beta_2)$, via \Ba-,~\Va-,
  and~\MPa-moves (leading to $(\calT'_1,\beta'_1)$ and
  $(\calT'_2,\beta'_2)$, respectively);

\item an application of the relative version of the Alexander theorem
  to relate $(\calT'_1,\beta'_1)$ and $(\calT'_2,\beta'_2)$ via
  \Ba-~and \MPa-moves;

\item an elimination of each \Ba-move by substituting it with a
  \Ca-move.

\end{itemize}
We first recall the relative version of the Alexander theorem, and
then we describe each of the three steps.

\subsection{Alexander's theorem}

As said above, our proof relies on the relative version of Alexander's
theorem, so we recall it (the proof is quite easy and can be found
in~\cite{Turaev-Viro}).
Let us consider a simplex $\sigma$ of a polyhedron $|\calP|$ with a
(non-loose) triangulation $\calP$.
Let us define a move on the triangulation $\calP$: the substitution of
the closed star, $\clst{\sigma}$, of $\sigma$ with the cone on
$\partial\clst{\sigma}$ with respect to a point in the interior of
$\sigma$ will be called A-{\em move}; the inverse of an A-move will be
also called an A-move.
Note that an A-move does not change the homeomorphism class of
$|\calP|$, but only the triangulation $\calP$.
The following theorem states that A-moves are enough to obtain all the
triangulations of $|\calP|$ from any given one, leaving fixed a
sub-polyhedron $|\calQ|$.
\begin{teo}\label{Alexander:teo}
Let $|\calP|$ be a dimensionally homogeneous polyhedron and let
$|\calQ|$ be a sub-polyhedron of $|\calP|$.
Then two triangulations of $|\calP|$, whose restrictions to $|\calQ|$
coincide, can be obtained from each other via a sequence of {\rm
  A}-moves which do not change the triangulation of $|\calQ|$.
\end{teo}

\paragraph{Reduction to B-~and MP-moves}
Now we prove a modification of a result due to Pachner (Theorem~4.14
of~\cite{Pachner}), that he stated only for manifolds.
Let us call {\em singular manifold with boundary} a finite polyhedron
$|\calP|$ such that the link of every point (of $|\calP|$) is a
surface with (possibly empty) boundary.
Such a space is the generalization with boundary of the so called {\em
  singular manifolds}.
In fact, we have an obvious definition of the {\em boundary} $\partial
|\calP|$ of $|\calP|$ as the 2-dimensional sub-polyhedron of $|\calP|$
made of the closure of the triangles lying in only one tetrahedron.
Obviously, in $|\calP|$ there are only a finite number of points
having link different from the 2-sphere or the 2-disk.
We have the following corollary of Theorem~\ref{Alexander:teo}.
\begin{prop}\label{gener_Pachner:prop}
Let $|\calP|$ be a singular manifold with boundary and let $|\calQ|$
be a sub-polyhedron of $|\calP|$ containing $\partial |\calP|$.
Then two triangulations of $|\calP|$, whose restrictions to $|\calQ|$
coincide, can be obtained from each other via a sequence of {\rm
  B}-~and {\rm MP}-moves, which do not change the common triangulation
of $|\calQ|$.
\end{prop}
\dimo{gener_Pachner:prop}
By Theorem~\ref{Alexander:teo}, the two triangulations can be obtained
from each other via a sequence of A-moves which do not change the
triangulation of $|\calQ|$.
To conclude the proof, we show that each A-move in this sequence is a
composition of B-~and MP-moves which do not change the triangulation
of $|\calQ|$.
There are four different types of A-move depending on the dimension of
the simplex $\sigma$ the A-move is applied to.
\begin{description}

\item{${\rm dim}(\sigma)=0$.}
This case is obvious; in fact, $\sigma$ is a vertex, so
$\clst{\sigma}$ is already the cone on $\partial\clst{\sigma}$ with
respect to $\sigma$, and the A-move is the identity.

\item{${\rm dim}(\sigma)=1$.}
Here $\sigma$ is an edge, so the A-move on $\sigma$ ``divides''
$\sigma$ adding a vertex as shown in Fig.~\ref{A_move_edge:fig}.
\begin{figure}
  \centerline{\includegraphics{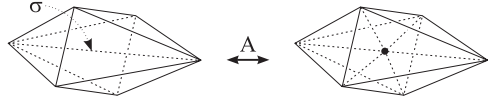}}
  \caption{The A-move on the edge $\sigma$ (with four tetrahedra in
    $\st{\sigma}$).}
  \label{A_move_edge:fig}
\end{figure}
Consider the open star, $\st{\sigma}$, of $\sigma$ shown in
Fig.~\ref{A_move_edge:fig}-left.
Note that $\st{\sigma}$ contains at least three tetrahedra: we
describe the case for four tetrahedra, other cases being similar.
The A-move is the composition of the moves shown in
Fig.~\ref{A_move_edge_moves:fig}: one positive B-move, two positive
MP-moves, and one negative MP-move.
\begin{figure}
  \centerline{\includegraphics{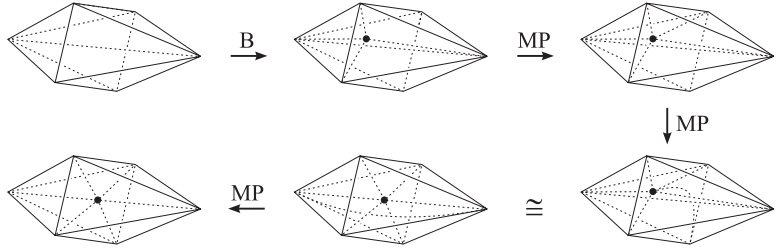}}
  \caption{The A-move on an edge is a composition of B-~and MP-moves.}
  \label{A_move_edge_moves:fig}
\end{figure}

\item{${\rm dim}(\sigma)=2$.} 
In Fig.~\ref{A_move_tria_moves:fig} we show that the A-move on a
triangle is a composition of one positive B-move and one positive
MP-move.
\begin{figure}
  \centerline{\includegraphics{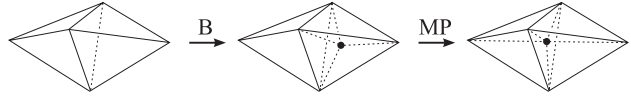}}
  \caption{The A-move on a triangle is a composition of B-~and
    MP-moves.}
  \label{A_move_tria_moves:fig}
\end{figure}

\item{${\rm dim}(\sigma)=3$.}
The A-move is already a B-move.

\end{description}
Finally, note that all the B-~and MP-moves described above do not
change the common triangulation of $|\calQ|$.
\finedimo{gener_Pachner:prop}

\subsection{Desingularization}\label{desingularization:subsec}

Let $(\calT,\beta)$ be an ideal triangulation of a pair $(M,\alpha)$.
As said above, the idea is to eliminate the singularities of the loose
triangulation $(\widehat{\calT},\widehat{\beta})$, via \Ba-,~\Va-, and
\MPa-moves, to be able to apply Proposition~\ref{gener_Pachner:prop}.
We will see that we cannot eliminate all the singularities, because we
cannot eliminate the edges belonging to $\widehat{\beta}$.
Since $\widehat{\calT}$ is a loose triangulation (of $\widehat{M}$),
there could be a singular edge of $\widehat{\calT}$ with coinciding
endpoints; such an edge will be called {\em loop}.
\begin{prop}\label{desingularization:prop}
Let $(\calT,\beta)$ be an ideal triangulation of a pair $(M,\alpha)$.
Then there exists an ideal triangulation $(\calT',\beta')$ of $(M
\setminus \cup B_k,\alpha)$, where the $B_k$'s are $3$-balls disjoint
from each other and from $\alpha$, such that the following facts
hold.
\newcounter{listi}
\begin{list}{{\rm \arabic{listi}.}}{\usecounter{listi} \setlength{\labelwidth}{1cm}}

\item $(\calT',\beta')$ is obtained from $(\calT,\beta)$ via
  \Ba-,~\Va-, and \MPa-moves.

\item The loose triangulation $(\widehat{\calT}',\widehat{\beta}')$
  has only the following types of singularities:
  \newcounter{listii}
  \begin{list}{{\rm (\alph{listii})}}{\usecounter{listii} \setlength{\labelwidth}{1cm}}
  \item an edge $(\widehat{\beta}')^{(i)}$ which is a loop,
  \item a pair of edges $(\widehat{\beta}')^{(i)}$ sharing both the
    endpoints,
  \item a pair of edges (giving a multiple adjacency) in
    $\clst{(\widehat{\beta}')^{(i)}}$ if $(\widehat{\beta}')^{(i)}$ is
    a loop.
  \end{list}

\item Each $(\widehat{\beta}')^{(i)}$ has a neighborhood
  $\calN((\widehat{\beta}')^{(i)})$ such that:
  \begin{list}{{\rm (\alph{listii})}}{\usecounter{listii} \setlength{\labelwidth}{1cm}}
  \item if $(\widehat{\beta}')^{(i)}$ is not a loop,
    $\calN((\widehat{\beta}')^{(i)})$ is made of exactly three
    tetrehedra;
  \item if $(\widehat{\beta}')^{(i)}$ is a loop,
    $\calN((\widehat{\beta}')^{(i)})$ is the cone on a triangle
    $\theta$, where $\theta$ is triangulated as shown in
    Fig.~\ref{triangulated_triangle:fig}, the endpoints of the cone on
    the barycentre $b$ of $\theta$ are identified together, and the
    loop $(\widehat{\beta}')^{(i)}$ is just this edge with identified
    endpoints.
    \begin{figure}
      \centerline{\includegraphics{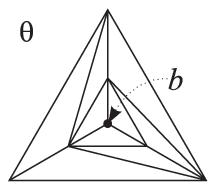}}
      \caption{The triangle $\theta$ and its barycentre $b$.}
      \label{triangulated_triangle:fig}
    \end{figure}
  \end{list}

\item $\calN((\widehat{\beta}')^{(i)}) \cap
  \calN((\widehat{\beta}')^{(j)}) = (\widehat{\beta}')^{(i)} \cap
  (\widehat{\beta}')^{(j)}$ for each $i \neq j$, and
  $\calN((\widehat{\beta}')^{(i)}) \cap \widehat{\partial M} =
  (\widehat{\beta}')^{(i)} \cap \widehat{\partial M}$ for $i=1,\ldots
  ,n$.

\end{list}
\end{prop}

\dimo{desingularization:prop}
The loose triangulation $(\widehat{\calT},\widehat{\beta})$ has
different types of singularity: we eliminate the singularities type by
type, being careful not to create any singularity of the types already
eliminated.
Note that we need to analyze only the singularities for tetrahedra,
because both a singular triangle and a singular edge are contained in
a singular tetrahedron.
There are 6 different types of singularity for tetrahedra.
For the sake of shortness, we continue calling $(\calT,\beta)$ also
the triangulations obtained during the proof, also if they are
actually different from $(\calT,\beta)$.

\begin{figure}
  \centerline{\includegraphics{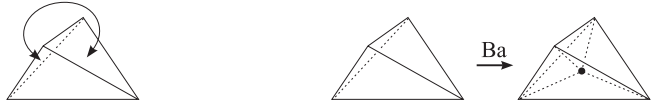}}
  \caption{A tetrahedron self-adjacent along triangles (left) and the
    \Ba-move which eliminates the self-adjacency (right).}
  \label{self_adj_tria:fig}
\end{figure}

\paragraph{Self-adjacency along triangles}
The tetrahedron is shown in Fig.~\ref{self_adj_tria:fig}-left; the
\Ba-move eliminating the self-adjacency is shown in
Fig.~\ref{self_adj_tria:fig}-right.
Note that no new self-adjacency along triangles has been created.

\paragraph{Self-adjacency along edges}
This case is more complicated than the previous one: for each
tetrahedron the number of edges which are identified together can vary
between 2 and 6.
An easy induction on the maximal number of edges identified together
in a tetrahedron and on the number of tetrahedra having such a maximal
number of identifications reduces the number of cases to two.
\begin{enumerate}

\item If two edges which are identified together do not share any
  vertex (in the unfolded version of the tetrahedron), a positive
  \Ba-move is enough to eliminate the singularity, see
  Fig.~\ref{self_adj_edges_1:fig}.
  \begin{figure}
    \centerline{\includegraphics{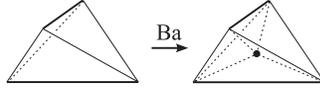}}
    \caption{Elimination of self-adjacency along edges (first case).
      The edges which are identified together are drawn thick.}
  \label{self_adj_edges_1:fig}
\end{figure}
 
\item If two edges which are identified together share a vertex (in
  the unfolded version of the tetrahedron) the situation is slightly
  more difficult.
  Let us start by calling $T$ the tetrahedron.
  Note that the tetrahedron $T'$, attached to $T$ along the triangle
  containing the two identified edges, is different from $T$, because we
  have already eliminated the self-adjacencies of tetrahedra along
  triangles.
  So a positive \Ba-move and a positive \MPa-move can be applied to
  eliminate the self-adjacency, see Fig.~\ref{self_adj_edges_2:fig}.
  \begin{figure}
    \centerline{\includegraphics{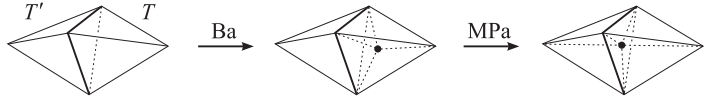}}
    \caption{Elimination of self-adjacency along edges (second case).
      The edges which are identified together are drawn thick.}
    \label{self_adj_edges_2:fig}
  \end{figure}
  
\end{enumerate}
Note that no new self-adjacency along either triangles or edges has
been created.

\paragraph{Multiple adjacency along triangles or edges}
The situation is analogous to the case of self-adjacencies along
triangles or edges, respectively; so it can be treated similarly.

\paragraph{}
Before continuing desingularization, we modify the loose triangulation
obtained after the first part of the process to ``isolate'' each edge
$\widehat{\beta}^{(i)}$ of $\widehat{\beta}$.
Namely, we apply \Ba-~and \MPa-moves to obtain point~4 of the
statement, {\em i.e.}~$\calN((\widehat{\beta}')^{(i)}) \cap
\calN((\widehat{\beta}')^{(j)}) = (\widehat{\beta}')^{(i)} \cap
(\widehat{\beta}')^{(j)}$ for each $i \neq j$, and
$\calN((\widehat{\beta}')^{(i)}) \cap \widehat{\partial M} =
(\widehat{\beta}')^{(i)} \cap \widehat{\partial M}$ for $i=1,\ldots
,n$.
The situation is similar to desingularization: we eliminate the
intersections between two $\clst{\widehat{\beta}^{(i)}}$'s and between
each $\clst{\widehat{\beta}^{(i)}}$ and $\widehat{\partial M}$ step by
step, being careful not to add any intersection of the types already
eliminated.
First we eliminate the intersections between each
$\clst{\widehat{\beta}^{(i)}}$ and $\widehat{\partial M}$.
If $\clst{\widehat{\beta}^{(i)}} \cap \widehat{\partial M}$ contains a
vertex $v$ different from the endpoints of the edge
$\widehat{\beta}^{(i)}$, then we perform the moves already described
to eliminate the self-adjacency of tetrahedra along edges (second
case); so $v$ belongs no more to $\clst{\widehat{\beta}^{(i)}} \cap
\widehat{\partial M}$.
Let us consider now the intersection between two
$\clst{\widehat{\beta}^{(i)}}$'s.
They may share (out of the intersection between the edges
$\widehat{\beta}^{(i)}$) tetrahedra, triangles, edges and vertices
(different from the endpoints of the edges $\widehat{\beta}^{(i)}$). 
\begin{description}
\item{{\em Tetrahedra}:} if two $\clst{\widehat{\beta}^{(i)}}$'s share
  a tetrahedron, we note that the two $\widehat{\beta}^{(i)}$'s belong
  to one tetrahedron, so we perform the moves already described to
  eliminate self-adjacency of tetrahedra along edges.

\item{{\em Triangles}:} if the common simplex is a triangle, we
  perform the move already used to eliminate multiple adjacency of
  tetrahedra along triangles.

\item{{\em Edges}:} if the common simplex is an edge, we perform the
  moves already used to eliminate multiple adjacency of tetrahedra
  along edges.

\item{{\em Vertices}:} if two $\clst{\widehat{\beta}^{(i)}}$'s share a
  vertex (different from the endpoints of the edges
  $\widehat{\beta}^{(i)}$), we perform the moves already described to
  eliminate the self-adjacency of tetrahedra along edges (second
  case).
\end{description}
Note that all the moves described above are admissible, and that no
new singularity of the types already eliminated has been created.
Let us continue now with desingularization.

\paragraph{Self-adjacency along vertices}
If two vertices of a tetrahedron are identified together, let us call
$e$ the edge which is a loop (if there is more than one edge like $e$,
we repeat the procedure).
There are two cases depending on whether the edge $e$ belongs to
$\widehat{\beta}$ or not.

\subparagraph{First case:~$e \not \hspace{-1pt} \in \widehat{\beta}$}
Consider the unfolded version of $\clst{e}$: the case for four
tetrahedra is shown in Fig.~\ref{A_move_edge:fig}-left.
We know that $\clst{e}$ contains at least three tetrahedra, because we
have already eliminated self-adjacencies and multiple adjacencies of
tetrahedra along triangles.
The idea is now to ``divide'' the edge $e$ by adding a vertex, as
shown in Fig.~\ref{A_move_edge:fig}.
The situation is analogous to that of the proof of
Proposition~\ref{gener_Pachner:prop} when the case of ${\rm
  dim}(\sigma)=1$ is analyzed; the only difference is that now some
boundary faces of $\clst{e}$ could be glued together, but this does
not matter: we can repeat the same B-~and MP-moves, ``dividing'' the
edge $e$, as shown in Fig.~\ref{A_move_edge_moves:fig}.
We conclude by noting that each move is admissible: the first three
are positive and the last one eliminates the edge $e$ which does not
belong to $\widehat{\beta}$.
Note also that no new singularity of the types already eliminated has
been created.

\subparagraph{Second case:~$e \in \widehat{\beta}$}
For the sake of clarity, let us call $\widehat{\beta}^{(i)}$ the edge
$e$.
Note that we cannot eliminate the singularity: in fact we cannot
eliminate the edge $\widehat{\beta}^{(i)}$, so each tetrahedron in
$\clst{\widehat{\beta}^{(i)}}$ always has $\widehat{\beta}^{(i)}$ as
an edge and it is always singular.
But we will modify a neighborhood of $\widehat{\beta}^{(i)}$ to obtain
point~3b of the statement.
Recall that $\clst{\widehat{\beta}^{(i)}} \setminus
(\widehat{\beta}^{(i)} \cap \widehat{\partial M})$ and
$\clst{\widehat{\beta}^{(j)}} \setminus (\widehat{\beta}^{(j)} \cap
\widehat{\partial M})$ are disjoint for each $j \neq i$.
First we will modify $\clst{\widehat{\beta}^{(i)}}$ via \Ba-,~\Va-,
and \MPa-moves to have that $\clst{\widehat{\beta}^{(i)}}$ is made of
exactly three tetrahedra; then we will modify these tetrahedra to
obtain point~3b of the statement.

\begin{figure}
  \centerline{\includegraphics{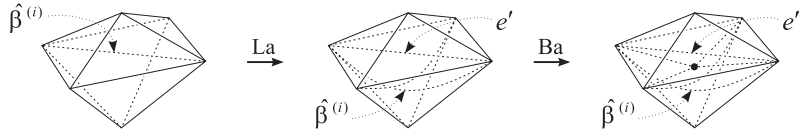}}
  \caption{How to simplify $\clst{\widehat{\beta}^{(i)}}$ so to have
    only three tetrahedra in it (case of four tetrahedra in
    $\clst{\widehat{\beta}^{(i)}}$).
    The endpoints of $\widehat{\beta}^{(i)}$ are identified together.}
  \label{link_edge_alpha:fig}
\end{figure}

Let us describe the first modification of
$\clst{\widehat{\beta}^{(i)}}$.
Note that $\clst{\widehat{\beta}^{(i)}}$ cannot be made of one or two
tetrahedra because we have already eliminated self-adjacencies and
multiple adjacencies of tetrahedra along triangles.
So let us suppose that $\clst{\widehat{\beta}^{(i)}}$ is made of at
least four tetrahedra and let us modify the loose triangulation to
have that $\clst{\widehat{\beta}^{(i)}}$ is made of three tetrahedra.
For the sake of clarity, in Fig.~\ref{link_edge_alpha:fig} we have
shown only the case of four tetrahedra in
$\clst{\widehat{\beta}^{(i)}}$: the other cases are analogous.
We apply a positive \La-move (which is a composition of \Va-~and
\MPa-moves), choosing to leave in $\widehat{\beta}$ the edge whose
star is made of three tetrahedra; we eliminate the multiple adjacency
created by the \La-move with a positive \Ba-move; we eliminate the
singularity of the edge $e'$ (``parallel'' to $\widehat{\beta}^{(i)}$)
created by the \La-move as we have done above ($e' \notin
\widehat{\beta}$).

\begin{figure}
  \centerline{\includegraphics{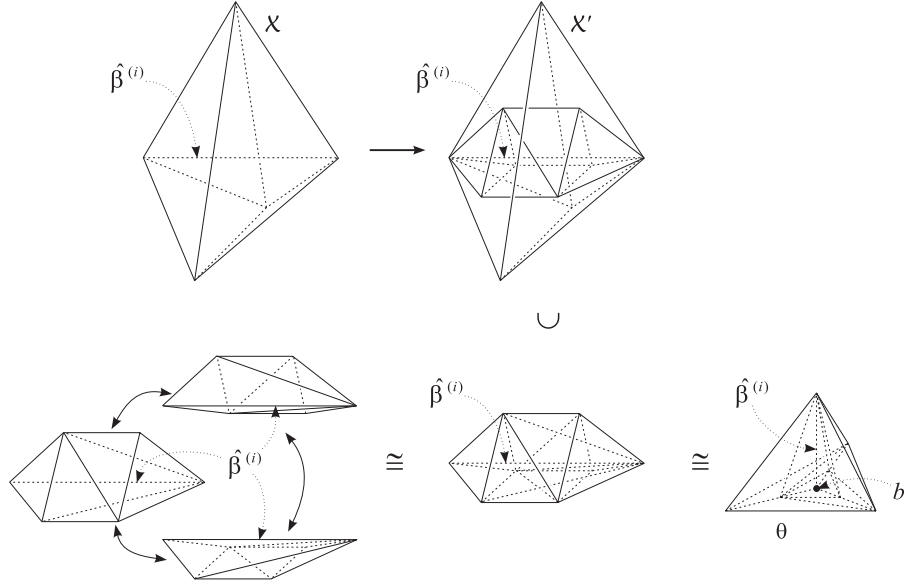}}
  \caption{The old $\clst{\widehat{\beta}^{(i)}}$ is modified so that
    the new $\calN(\widehat{\beta}^{(i)})$ is contained in the old
    $\clst{\widehat{\beta}^{(i)}}$ (shown transparent).
    We do not show the whole of $\calX'$: we show only how it appears
    near $\widehat{\beta}^{(i)}$.
    The endpoints of $\widehat{\beta}^{(i)}$ are identified together.}
  \label{star_edge_self_alpha:fig}
\end{figure}

Let us pass to the second modification of
$\clst{\widehat{\beta}^{(i)}}$, which is now made of three
tetrahedra.
Consider the unfolded version of $\clst{\widehat{\beta}^{(i)}}$: it
can be seen as a triangulation, say $\calX$, of the 3-ball, see
Fig.~\ref{star_edge_self_alpha:fig}-left.
Let $\calX'$ be another triangulation of the 3-ball such that:
\begin{itemize}

\item $\calX$ and $\calX'$ coincide on the boundary of the 3-ball and
  on the edge $\widehat{\beta}^{(i)}$;

\item $\calX'$ appears, near $\widehat{\beta}^{(i)}$, as in
  Fig.~\ref{star_edge_self_alpha:fig}-right;

\item any two boundary faces of $\calX'$ do not belong to the same
  tetrahedron.

\end{itemize}
It is very easy to find such an $\calX'$.
Now, $\calX$ and $\calX'$ have in common the boundary and the edge
$\widehat{\beta}^{(i)}$, so we can apply
Proposition~\ref{gener_Pachner:prop} to obtain $\calX'$ from $\calX$
via B-~and MP-moves not involving both the edge
$\widehat{\beta}^{(i)}$ and the boundary.
Repeating these moves on the folded version of $\calX$ contained in
$\calT$, we substitute it with a folded version of $\calX'$ using
B-~and MP-moves which are admissible because they have support in the
folded version of $\calX$ and do not involve the edge
$\widehat{\beta}^{(i)}$.
Now a neighborhood of $\widehat{\beta}^{(i)}$, say
$\calN(\widehat{\beta}^{(i)})$ appears as in
Fig.~\ref{star_edge_self_alpha:fig}-bottom.
Note that $\calN(\widehat{\beta}^{(i)})$ is the cone on the triangle
$\theta$ shown in Fig.~\ref{triangulated_triangle:fig} where the
endpoints of the cone on the barycentre $b$ are identified together,
that $(\widehat{\beta})^{(i)}$ is just this edge with identified
endpoints, and that no new singularity of the types already eliminated
has been created.

\paragraph{Multiple adjacency along vertices}
The situation is analogous to the case of self-adjacency along
vertices, but there are some differences to point out.
The idea is to ``divide'' one of the edges giving the singularity, so
the moves to apply are those applied to eliminate self-adjacency along
vertices when $e \notin \widehat{\beta}$.
But there are two exceptions.
\begin{enumerate}
\item We cannot ``divide'' the edges belonging to $\widehat{\beta}$ so
  we cannot eliminate the singularity created by two edges of
  $\widehat{\beta}$ sharing both the endpoints.

\item If an edge $(\widehat{\beta})^{(i)}$ (belonging to
  $\widehat{\beta}$) is a loop, then we do not divide any of the edges
  belonging to the closed star of $\widehat{\beta}^{(i)}$, because
  such an edge has in its closed star a loop (the edge
  $\widehat{\beta}^{(i)}$) and the moves described above would create
  a new multiple adjacency.
\end{enumerate}
For the other cases we can eliminate the multiple adjacency as we have
done for self-adjacencies along vertices with $e \notin
\widehat{\beta}$, because both the moves are admissible and we do not
add any of the singularities of the types already eliminated.
Finally, let us deal with the two exceptions.
\begin{enumerate}
\item For each edge $\widehat{\beta}^{(i)}$ which is not a loop, we
  modify $\clst{\widehat{\beta}^{(i)}}$ to have that it is made of
  exactly three tetrahedra, as we have done above for the first
  modification of $\clst{\widehat{\beta}^{(i)}}$ for the
  $\widehat{\beta}^{(i)}$'s which are loops.

\item We do not operate on the edges belonging to the closed star of
  the $\widehat{\beta}^{(i)}$'s which are loops.
\end{enumerate}

\paragraph{Conclusion}
Repeating the moves described above on the ideal triangulation $(\calT
,\beta)$ of the pair $(M,\alpha)$, we obtain, via \Ba-,~\Va-, and
\MPa-moves, an ideal triangulation $(\calT',\beta')$ of $(M \setminus
\cup B_k, \alpha)$, where the $B_k$'s are $3$-balls disjoint from each
other and from $\alpha$.
We have eliminated almost all the singularities of $\widehat{\calT}$,
but there are three types of singularity we cannot eliminate (those
due to $\alpha$).
These three types of singularity are exactly those described in
point~2 of the statement.
The check that $(\calT',\beta')$ is the desired ideal triangulation is
straight-forward, so we leave it to the reader.
\finedimo{desingularization:prop}

\subsection{Application of the Alexander theorem}

Let us state (and prove) now a first result, which is a weak version
of Theorem~\ref{gener_MP:teo}.
\begin{prop}\label{gener_MP_weak:prop}
Two marked ideal triangulations of a pair $(M,\alpha)$ can be obtained
from each other via a sequence of \Ba-,~\Va-, and \MPa-moves, such
that the negative \Ba-moves do not eliminate the spherical boundary
components of $\partial M$.
\end{prop}
\dimo{gener_MP_weak:prop}
Let $(\calT_1,\beta_1)$ and $(\calT_2,\beta_2)$ be two ideal
triangulations of $(M,\alpha)$.
Let us apply Proposition~\ref{desingularization:prop} to each
$(\calT_i,\beta_i)$ obtaining $(\calT'_i,\beta'_i)$.
Recall that each $(\calT'_i,\beta'_i)$ is obtained from the
corresponding $(\calT_i,\beta_i)$ via \Ba-,~\Va-, and \MPa-moves, that
each $(\beta'_i)^{(j)}$ has a particular neighborhood
$\calN((\beta'_i)^{(j)})$, and that the loose triangulations
$(\widehat{\calT}'_i,\widehat{\beta}'_i)$ are almost desingularized
(the singularities are contained in the open neighborhood
$\inter{\calN((\widehat{\beta}'_i)^{(j)})}$).
Moreover, recall that the \Ba-move does not involve the spherical
boundary components of $\partial M$.
Obviously, since we have $\calN((\widehat{\beta}'_i)^{(j)}) \cap
\calN((\widehat{\beta}'_i)^{(k)}) = (\widehat{\beta}'_i)^{(j)} \cap
(\widehat{\beta}'_i)^{(k)}$ for each $j \neq k$, and
$\calN((\widehat{\beta}'_i)^{(j)}) \cap \widehat{\partial M} =
(\widehat{\beta}'_i)^{(j)} \cap \widehat{\partial M}$ for $j=1,\ldots
,n$, we can suppose, up to isotopy, that $\calN(\widehat{\beta}'_1)$
and $\calN(\widehat{\beta}'_2)$ coincide.

The strategy will now be to prove that $(\calT'_1,\beta'_1)$ and
$(\calT'_2,\beta'_2)$ are obtained from each other via \Ba-~and
\MPa-moves.
To do this, we will apply Proposition~\ref{gener_Pachner:prop} to
$\widehat{M} \setminus (\sqcup_j
\inter{\calN((\widehat{\beta}'_1)^{(j)})}) = \widehat{M} \setminus
(\sqcup_j \inter{\calN((\widehat{\beta}'_2)^{(j)})})$.
Since the singularities of the loose triangulations
$\widehat{\calT}'_i$ are contained in the
$\inter{\calN((\widehat{\beta}'_i)^{(j)})}$'s (see
Proposition~\ref{desingularization:prop}), the triangulations
$\widehat{\calT}'_i \setminus (\sqcup_j
\inter{\calN((\widehat{\beta}'_i)^{(j)})})$ are actually non-loose.
Moreover, the two $\widehat{\calT}'_i \setminus (\sqcup_j
\inter{\calN((\widehat{\beta}'_i)^{(j)})})$'s coincide on the boundary
and on $\widehat{\partial M}$.
Then, we can apply Proposition~\ref{gener_Pachner:prop} to transform
$\widehat{\calT}'_1 \setminus (\sqcup_j
\inter{\calN((\widehat{\beta}'_1)^{(j)})})$ into $\widehat{\calT}'_2
\setminus (\sqcup_j \inter{\calN((\widehat{\beta}'_2)^{(j)})})$ via
B-~and MP-moves having support out of $\widehat{\partial M}$.
Obviously, these moves can be applied on the loose triangulation
$\widehat{\calT}'_1$ transforming it into $\widehat{\calT}'_2$, they
are all admissible, and they transform the loose triangulation
$(\widehat{\calT}'_1,\widehat{\beta}'_1)$ into
$(\widehat{\calT}'_2,\widehat{\beta}'_2)$; moreover, the negative
\Ba-moves do not eliminate the points belonging to $\widehat{\partial
  M}$.
The desired sequence is obtained by repeating the moves on the ideal
triangulations $(\calT'_i,\beta'_i)$ of $(M,\alpha)$.
\finedimo{gener_MP_weak:prop}

\subsection{Elimination of Ba-moves}

To deduce Theorem~\ref{gener_MP:teo} from
Proposition~\ref{gener_MP_weak:prop}, we generalize an idea of
Matveev~\cite{Matveev:calculus} to the setting of marked spines.

\dimo{gener_MP:teo}
Let $(\calT_1,\beta_1)$ and $(\calT_2,\beta_2)$ be two ideal
triangulations of $(M,\alpha)$.
By Proposition~\ref{gener_MP_weak:prop}, we have that
$(\calT_2,\beta_2)$ is obtained from $(\calT_1,\beta_1)$ via
\Ba-,~\Va-, and \MPa-moves, such that the negative \Ba-moves do not
eliminate the spherical boundary components of $\partial M$.
The idea of the proof consists in replacing each \Ba-move with a
\Ca-move, and each \Va-~or \MPa-move with suitable sequences of
\La-,~\Va-, and \MPa-moves.
Let us pass to the dual spine viewpoint: for $i=1,2$, let
$(P_i,\beta_i)$ be the spine dual to $(\calT_i,\beta_i)$.

First of all, note that in the passages along the sequence of
\Ba-,~\Va-,~and \MPa-moves we get (standard) spines $P_*$ of $M$ minus
some balls; so each $M \setminus P_*$ is a disjoint union of $\partial
M \times [0,1)$ and some balls.
When a positive \Ba-move is applied, a proper ball $\calB$ appears.
Let us continue calling {\em proper ball} (and continue indicating it
by $\calB$) its transformations after the others \Ba-,~\Va-, and
\MPa-moves, until it disappears because of a negative \Ba-move (each
proper ball has to disappear).
Note that, conversely, the negative \Ba-moves eliminate only the
proper balls.
Note also that each $\calB$ is an open ball with boundary contained in
$P_*$ and it is not touched by the edges belonging to $\alpha$.

We will not replace all the \Ba-moves (with \Ca-moves) at the same
time; instead, we will concentrate on one positive \Ba-move and on the
negative \Ba-move eliminating the proper ball created by the positive
\Ba-move.
The strategy will be to replace these two \Ba-moves with two
\Ca-moves, any other \Ba-move with a suitable sequence of only one
\Ba-move and \La-,~\Va-, and \MPa-moves, and each \Va-~or \MPa-move
with a suitable sequence of \La-,~\Va-, and \MPa-moves only.
In such a way we will decrease, by two, the number of \Ba-moves in the
sequence.
By repeating this procedure we can eliminate all the \Ba-moves and we
can complete the proof.

Let us describe the procedure in details.
Let $s$ be the following sequence of \Ba-,~\Va-, and \MPa-moves
transforming $(P_1,\beta_1)$ into $(P_2,\beta_2)$:
\begin{eqnarray*}
(P_1,\beta_1)
\stackrel{s_1}{\longrightarrow}
(Q_0,\eta_0)
\stackrel{\Ba^+}{\longrightarrow}
(Q_1,\eta_1)
\stackrel{m_1}{\longrightarrow}
(Q_2,\eta_2)
\stackrel{m_2}{\longrightarrow}
\quad\ldots \\
\ldots\quad
\stackrel{m_{r-1}}{\longrightarrow}
(Q_r,\eta_r)
\stackrel{\Ba^-}{\longrightarrow}
(Q_{r+1},\eta_{r+1})
\stackrel{s_2}{\longrightarrow}
(P_2,\beta_2),
\end{eqnarray*}
where $s_1$ and $s_2$ are sequences of moves we will not replace,
$\Ba^+$ (respectively, $\Ba^-$) is the positive (respectively,
negative) move we will replace with a positive (respectively,
negative) \Ca-move, and the $m_j$'s are the other moves we will
replace.
From now on, we will denote by $\calB$ both the proper ball created by
$\Ba^+$ (and eliminated by $\Ba^-$) and its transformations after the
$m_j$'s.
To decrease by two the number of \Ba-moves, we will find a sequence
$s'$ transforming $(P_1,\beta_1)$ into $(P_2,\beta_2)$ and appearing
as follows:
\begin{eqnarray*}
(P_1,\beta_1)
\stackrel{s_1}{\longrightarrow}
(Q_0,\eta_0)
\stackrel{\Ca^+}{\longrightarrow}
(\widetilde{Q}_1,\widetilde{\eta}_1)
\stackrel{\widetilde{m}_1}{\longrightarrow}
(\widetilde{Q}_2,\widetilde{\eta}_2)
\stackrel{\widetilde{m}_2}{\longrightarrow}
\quad\ldots \\
\ldots\quad
\stackrel{\widetilde{m}_{r-1}}{\longrightarrow}
(\widetilde{Q}_r,\widetilde{\eta}_r)
\stackrel{\Ca^-}{\longrightarrow}
(Q_{r+1},\eta_{r+1})
\stackrel{s_2}{\longrightarrow}
(P_2,\beta_2),
\end{eqnarray*}
where $s_1$ and $s_2$ are the same sequences as above, $\Ca^+$
(respectively, $\Ca^-$) is the positive (respectively, negative) move
replacing $\Ba^+$ (respectively, $\Ba^-$), and the $\widetilde{m}_j$'s
are sequences of moves (composed either by only one \Ba-move and some
\La-,~\Va-, and \MPa-moves if $m_j$ is a \Ba-move, or by only
\La-,~\Va-, and \MPa-moves otherwise) replacing the $m_j$'s.

Let us start by replacing $\Ba^+$ with a positive \Ca-move $\Ca^+$
(the position of the arch can be random).
After applying $\Ca^+$ to $(Q_0,\eta_0)$ we obtain a spine
$(\widetilde{Q}_1,\widetilde{\eta}_1)$ which differs from
$(Q_1,\eta_1)$ only for the presence of an arch connecting the proper
ball $\calB$ to $M \setminus (Q_1 \cup \calB)$, see
Fig.~\ref{C_move:fig}-right.
Note that the arch joins a region $R_1$ of $\partial \calB$ with
another one, $R_2$, of $Q_1$; if $R_2$ belongs to $\eta_1$, then $R_1$
is a part of $\partial \calB$ belonging to $\widetilde{\eta}_1$.
Note also that $R_1$ is the only part of $\partial \calB$ which can
belong to $\widetilde{\eta}_1$, and that $R_0$ does not belong to
$\widetilde{\eta}_1$.
Now the sequence $s'$ appears as follows:
\begin{eqnarray*}
(P_1,\beta_1)
\stackrel{s_1}{\longrightarrow}
(Q_0,\eta_0)
\stackrel{\Ca^+}{\longrightarrow}
(\widetilde{Q}_1,\widetilde{\eta}_1).
\end{eqnarray*}

The aim is now to replace the moves $m_j$.
If we try to apply $m_1$ also on
$(\widetilde{Q}_1,\widetilde{\eta}_1)$, we could fail because of the
presence of the arch created by the move $\Ca^+$.
So the idea is either to apply the move $m_j$ if the arch is not
involved in the move, or to move the arch before applying the move
otherwise.
To do this, we will use a recursive procedure.
Let $(Q_j,\eta_j)$ be a spine (of the sequence $s$) of $(M,\alpha)$
minus some balls (let us call $k$ the number of such balls).
Let $\calB$ be the connected component of $M \setminus Q_j$ containing
one of such balls.
Note that $\calB$ is an open ball embedded in $\inter{M}$, but its
closure $\overline{\calB}$ may not be a closed ball embedded in
$\inter{M}$.
In our recursive procedure $\calB$ is the proper ball created by the
move $\Ba^+$ and modified by the moves $m_i$, with $i<j$.
Let $(\widetilde{Q}_j,\widetilde{\eta}_j)$ be a spine of $(M,\alpha)$
minus $k-1$ balls, which differs from $(Q_j,\eta_j)$ only for the
presence of an arch connecting the proper ball $\calB$ to another
connected component of $M \setminus Q_j$.
Let moreover $m_j$ be an admissible move from $(Q_j,\eta_j)$ to
$(Q_{j+1},\eta_{j+1})$, which does not eliminate the proper ball
$\calB$.
Note that $(Q_{j+1},\eta_{j+1})$ is a spine of $(M,\alpha)$ minus $h$
balls, where $h=k-1,k,k+1$ depending on $m_j$.
Let us continue calling $\calB$ the transformation of $\calB$ under $m_j$.

The recursive pass consists in describing a sequence $\widetilde{m}_j$
of admissible moves (composed either by only one \Ba-move and some
\La-,~\Va-, and \MPa-moves if $m_j$ is a \Ba-move, or by only
\La-,~\Va-, and \MPa-moves otherwise) from
$(\widetilde{Q}_j,\widetilde{\eta}_j)$ to
$(\widetilde{Q}_{j+1},\widetilde{\eta}_{j+1})$, where
$(\widetilde{Q}_{j+1},\widetilde{\eta}_{j+1})$ is a spine of
$(M,\alpha)$ minus $h-1$ ball, which differs from
$(Q_{j+1},\eta_{j+1})$ only for the presence of an arch connecting the
ball $\calB$ to another connected component of $M \setminus Q_{j+1}$.
If $m_j$ can be applied ({\em i.e.}~the arch is far from the support
of $m_j$), then we apply $m_j$ to
$(\widetilde{Q}_j,\widetilde{\eta}_j)$ obtaining
$(\widetilde{Q}_{j+1},\widetilde{\eta}_{j+1})$, which obviously has
all the properties described above.
Note that there are some types of moves which can always be applied
because the arch is never involved, up to isotopy, in the move: such
moves are the positive \Ba-moves and the positive \Va-moves.
To replace the other \Ba-,~\Va-, and \MPa-moves, maybe we need to move
the arch so to be able to apply the move.
If $m_j$ cannot be applied (because of the presence of the arch), then
we move the arch before applying $m_j$.
Let us describe how to move the arch; afterwards we will continue the
substitution of $m_j$ with $\widetilde{m}_j$.

\paragraph{Arch-move}
Let $(\widetilde{Q}_j,\widetilde{\eta}_j)$ be the spine of
$(M,\alpha)$ minus $k-1$ balls, which has an arch we want to move.
Recall that $\calB$ is the proper ball connected to another connected
component of $M \setminus Q_j$ by the arch.
Moreover recall that the proper ball $\calB$ is an open ball embedded
in $M$, but (because of the moves $m_i$ with $i<j$) its closure
$\overline{\calB}$ could be not a closed ball embedded in $M$.

\begin{figure}
  \centerline{\includegraphics{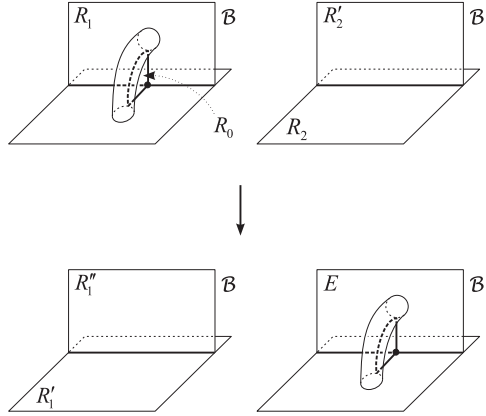}}
  \caption{The arch-move.
    We show on the left the situation near the arch we want to remove
    and on the right the situation near the point where we want to
    place the arch.}
  \label{arch_move:fig}
\end{figure}

We now define a spine $(\widetilde{Q}'_j,\widetilde{\eta}'_j)$ of
$(M,\alpha)$ minus $k-1$ balls.
Let $\widetilde{Q}'_j$ be the spine obtained from $\widetilde{Q}_j$ by
taking away the arch we want to move and by placing it in another
point, so that $\widetilde{Q}'_j$ is again a spine of $M$ minus $k-1$
balls and the ball $\calB$ is connected by the new arch to another
connected component of $M \setminus Q_j$, see
Fig.~\ref{arch_move:fig}.
The two conditions on $\widetilde{Q}'_j$ imply that the arch, after
the move, should be placed ``near'' $\partial \overline{\calB}$.
To define $\widetilde{\eta}'_j$, let us analyze the regions affected
by the move.
The region $R_1$ of $\widetilde{Q}_j$ (intersecting $\partial \calB$)
is divided (in $\widetilde{Q}'_j$) in two regions, $R'_1$ and
$R''_1$.
Note that these two regions belong also to the spine $(Q_j,\eta_j)$
and that only $R'_1$ can belong to $\eta_j$; if it belongs to
$\eta_j$, we impose to leave itself in $\widetilde{\eta}'_j$.
The little region $R_0$, which is eliminated by the arch-move, does
not belong to $\widetilde{\eta}_j$.
The other regions which are modified are the regions $R_2$ and $R'_2$,
which unite.
Note that these two regions belong also to $(Q_j,\eta_j)$ and that
only $R_2$ can belong to $\eta_j$ (because $R''_1$ is contained in
$\partial \calB$); if $R_2$ belongs to $\eta_j$, we impose to leave
the region $E$ in $\widetilde{\eta}'_j$.
The other regions are not modified, so we leave in
$\widetilde{\eta}'_j$ those belonging to $\widetilde{\eta}_j$.
Finally, note that $(\widetilde{Q}'_j,\widetilde{\eta}'_j)$ differs
from $(Q_j,\eta_j)$ only for the presence of the arch (connecting the
proper ball $\calB$ to another connected component of $M \setminus
Q_j$).
The transformation of $(\widetilde{Q}_j,\widetilde{\eta}_j)$ into
$(\widetilde{Q}'_j,\widetilde{\eta}'_j)$ will be called {\em
  arch-move}.

\begin{figure}
  \centerline{\includegraphics{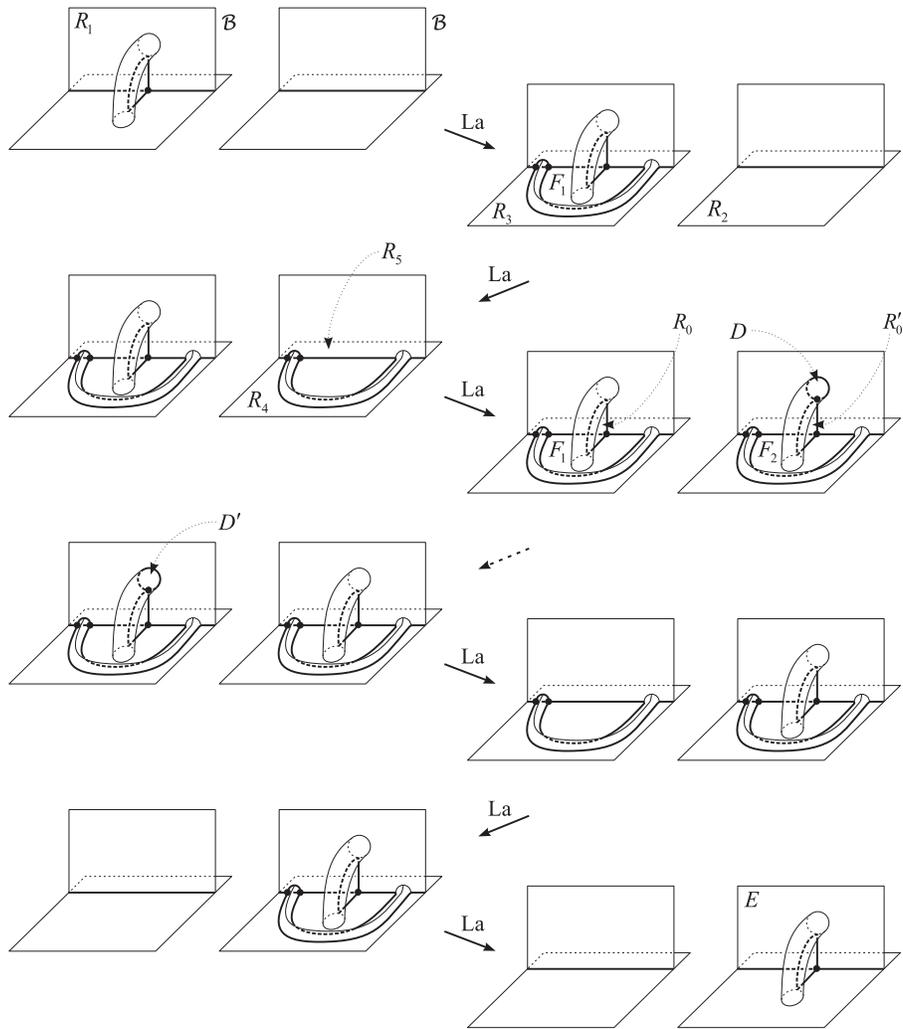}}
  \caption{The arch-move is a composition of \La-~and \MPa-moves.
    In each step, we show on the left the situation near the arch we
    want to remove and on the right the situation near the point where
    we want to place the arch.
    The dashed arrow denotes an isotopy of the little disc $D$.}
  \label{arch_L_MP:fig}
\end{figure}

Now we prove that each arch-move is a composition of \La-~and
\MPa-moves.
In Fig.~\ref{arch_L_MP:fig} we have shown the \La-~and \MPa-moves
transforming $(\widetilde{Q}_j,\widetilde{\eta}_j)$ into
$(\widetilde{Q}'_j,\widetilde{\eta}'_j)$: let us describe these
moves.
Note that the only region which can both intersect $\partial \calB$
and belong to $\widetilde{\eta}_j$ is $R_1$.
For the first positive L-move, if $R_1$ belongs to $\alpha$, we choose
to leave $R_3$ in $\alpha$.
Note that now no region in $\alpha$ intersects $\partial \calB$.
For the second positive L-move, if $R_2$ belongs to $\alpha$, we
choose to leave $R_4$ in $\alpha$.
The region $R_5$ does not belong to $\alpha$, because it intersects
$\partial \calB$, so the third positive L-move is admissible.
Let us now describe the move indicated by a dashed arrow.
Note that the proper ball $\calB$ can be seen as a tube $D^2 \times
[0,1]$, where $D^2 \times \{0\} = D$ and $D^2 \times \{1\} = D'$.
Obviously, we can move the disc $D$ through the tube from $D^2 \times
\{0\}$ to $D^2 \times \{1\}$ via an isotopy.
The move indicated by the dashed arrow consists exactly of this
isotopy of the little disc $D$ through the proper ball $\calB$: more
precisely, if one of the two arches (or both of them) are inside
$\calB$ (namely, the little discs $R_0$ and $R'_0$ are contained in
$\calB$), the isotopy is through $\calB$ minus both the arch and the
tube inside it.
At the end of the isotopy the little disc $D$ coincides with $D'$, so
it lays near the arch we want to remove.
A simple general position argument tells us that the isotopy can be
substituted with L-~and MP-moves, see Lemma~1.2.16
of~\cite{Matveev:new:book} for a precise proof.
All these moves are admissible because $F_1$, $F_2$, and the regions
intersecting $\partial \calB$ do not belong to $\alpha$.
For the same reason (and since $D'$ does not belong to $\alpha$) the
last three negative L-moves are admissible (obviously, the regions
united in each move are different).
To conclude, we note that the position of the regions in $\alpha$
after these moves is the same as after the arch-move.

\paragraph{Continuing substitution}
Recall that we want to replace the move $m_j$ which cannot be applied
to $(\widetilde{Q}_j,\widetilde{\eta}_j)$, because of the presence of
the arch.
We apply first an arch-move to $(\widetilde{Q}_j,\widetilde{\eta}_j)$
obtaining $(\widetilde{Q}'_j,\widetilde{\eta}'_j)$ and then the move
$m_j$.
Let us call $(\widetilde{Q}_{j+1},\widetilde{\eta}_{j+1})$ the spine
just obtained.
Note that, to apply the arch-move, we need to find a place where
placing the arch; but it is very easy to find such a place near
$\partial \overline{\calB}$ and far from the move $m_j$.
Note also that, by construction,
$(\widetilde{Q}_{j+1},\widetilde{\eta}_{j+1})$ differs from
$(Q_{j+1},\eta_{j+1})$ only for the presence of the arch (connecting
the proper ball $\calB $ to another connected component of $M
\setminus Q_{j+1}$).

With these substitutions, we have extended the sequence $s'$
obtaining:
\begin{eqnarray*}
(P_1,\beta_1)
\stackrel{s_1}{\longrightarrow}
(Q_0,\eta_0)
\stackrel{\Ca^+}{\longrightarrow}
(\widetilde{Q}_1,\widetilde{\eta}_1)
\stackrel{\widetilde{m}_1}{\longrightarrow}
(\widetilde{Q}_2,\widetilde{\eta}_2)
\stackrel{\widetilde{m}_2}{\longrightarrow}
\quad\ldots \\
\ldots\quad
\stackrel{\widetilde{m}_{r-1}}{\longrightarrow}
(\widetilde{Q}_r,\widetilde{\eta}_r).
\end{eqnarray*}

Let us consider now the move $\Ba^-$.
We have noted above that the spine
$(\widetilde{Q}_r,\widetilde{\eta}_r)$ differs from $(Q_r,\eta_r)$
only for the presence of the arch (connecting the proper ball $\calB$
to another connected component of $M \setminus Q_r$), so
$\widetilde{Q}_r$ near $\calB$ appears exactly as in
Fig.~\ref{C_move:fig}-centre.
Moreover, $R_1$ is the only part of $\partial\calB$ which can belong
to $\widetilde{\eta}_r$ and $R_0$ does not belong to
$\widetilde{\eta}_r$.
Obviously, a negative \Ca-move (which we call $\Ca^-$) can be applied
and the result is just $(Q_{r+1},\eta_{r+1})$.
Now the sequence $s'$ appears as follows:
\begin{eqnarray*}
(P_1,\beta_1)
\stackrel{s_1}{\longrightarrow}
(Q_0,\eta_0)
\stackrel{\Ca^+}{\longrightarrow}
(\widetilde{Q}_1,\widetilde{\eta}_1)
\stackrel{\widetilde{m}_1}{\longrightarrow}
(\widetilde{Q}_2,\widetilde{\eta}_2)
\stackrel{\widetilde{m}_2}{\longrightarrow}
\quad\ldots \\
\ldots\quad
\stackrel{\widetilde{m}_{r-1}}{\longrightarrow}
(\widetilde{Q}_r,\widetilde{\eta}_r)
\stackrel{\Ca^-}{\longrightarrow}
(Q_{r+1},\eta_{r+1}).
\end{eqnarray*}

To obtain the desired sequence, it is enough to complete the sequence
just obtained by composing it with the sequence $s_2$.
This proves the theorem.
\finedimo{gener_MP:teo}

\subsection{Another proof}\label{Bas_Ben_proof:subsec}

In this subsection we describe how Basehilac and Benedetti have
deduced Theorem~\ref{gener_MP:teo} from a result (due to Turaev and
Viro) which relies on the Matveev-Piergallini theorem.
For the sake of clarity, we describe the ideas of the proof, instead
of only stating Theorem~3.4.B of~\cite{Turaev-Viro}.
We restrict ourselves only to a sketch of the proof of
Theorem~\ref{gener_MP:teo}.

\vspace{1pt}\noindent{\it Sketch of the proof of} {\hspace{2pt}}\ref{gener_MP:teo}.
For $i=1,2$, let $(P_i,\beta_i)$ be the spine dual to an ideal
triangulation $(\calT_i,\beta_i)$.
Let $N(\alpha)$ be a little open regular neighborhood of $\alpha$ and
$M_\alpha = M \setminus N(\alpha)$.
Note that, up to choosing $N(\alpha)$ small with respect to $P_1$ and
$P_2$, we can suppose that $N(\alpha) \cap P_i = \cup_{j=1}^n
D_i^{(j)}$, where $D_i^{(j)}$ is an open disc with closure contained
in the (open) region $\beta_i^{(j)}$, for $i=1,2$ and $j=1,\ldots,n$.
Now, for $i=1,2$, we define two new polyhedra $Q_i$ and $R_i$ with
$Q_i \subset R_i \subset P_i$.
To get $Q_i$, we remove from $P_i$ all the (open) regions
$\beta_i^{(j)}$, and, to get $R_i$, we remove from $P_i$ all the open
discs $D_i^{(j)}$.
Note that we have a retraction $\pi_i$ of $M_\alpha$ onto $Q_i$.
Moreover, we have on $\partial M_\alpha$ a family $\lambda_i =
\{\lambda_i^{(1)},\ldots ,\lambda_i^{(n)}\}$ of disjoint simple
circles such that $\lambda_i^{(j)} = \partial D_i^{(j)} \subset
\beta_i^{(j)}$ and, up to isotopy, $R_i \setminus Q_i$ consists
precisely of the ``half-open'' annuli $\lambda_i^{(j)} \times [0,1)$
obtained by projecting $\lambda_i^{(j)}$ to $Q_i$ along $\pi_i$.
We have already described the ``inverse'' construction in
Subsection~\ref{id_tria:subsection} when we have proved existence of
marked ideal triangulations.
Of course any move on $R_i$ not affecting the $\lambda_i^{(j)}$'s
readily translates into an admissible move on $(P_i,\beta_i)$, and
conversely.
Obviously, up to isotopy, we can suppose that each $\lambda_1^{(j)}$
coincides with $\lambda_2^{(j)}$ and that each $D_1^{(j)}$ coincides
with $D_2^{(j)}$: let us call simply $\lambda^{(j)}$ the curve
$\lambda_1^{(j)}=\lambda_2^{(j)}$, $\lambda$ the collection
$\{\lambda^{(1)},\ldots,\lambda^{(n)}\}$, and $D^{(j)}$ the disc
$D_1^{(j)}=D_2^{(j)}$.

Now, $Q_i$ needs not to be standard, but one readily sees that
standardness can be achieved using C-~and L-moves on $R_i$ not
affecting the $\lambda^{(j)}$'s.
Now, $Q_1$ and $Q_2$ are standard spines of $M \setminus N(\alpha)$,
so, by Matveev-Piergallini theorem (see
Theorem~\ref{MP_calculus:teo}), we can transform $Q_1$ into $Q_2$ via
a deformation $Q_t$ (with $t \in [1,2]$) with elementary accidents
which are L-~and MP--moves.
Obviously, we can suppose that the elementary accidents occur at
different times.
Note that the $Q_t$'s are all quasi-standard spines, except for a
finite number of times when elementary accidents occur so
quasi-standardness is lost.

Parallelly, we have a deformation $\pi_t$ of $\pi_1$ into $\pi_2$,
where each $\pi_t$ is a retraction of $M \setminus N(\alpha)$ onto
$Q_t$.
Obviously, the annuli $[\lambda^{(j)},\pi_1(\lambda^{(j)}))$ are
transformed into $[\lambda^{(j)},\pi_2(\lambda^{(j)}))$ via annuli
$[\lambda^{(j)},\pi_t(\lambda^{(j)}))$.
By a general position argument, we can suppose that the accidents
occurring to $[\lambda,\pi_t(\lambda)) \cup Q_t$ are L-,~MP-, and
false L-moves not affecting the $\lambda^{(j)}$'s, where a {\em false}
L-move is a negative L-move not preserving standardness (actually it
is not an L-move).

Now, we have obtained a sequence of L-,~MP-, and false L-moves not
affecting the $\lambda^{(j)}$'s transforming $R_1$ into $R_2$.
To eliminate the false L-moves, we can use the same technique used in
Theorem~1.2.30 of~\cite{Matveev:new:book}, which states the following
(we use our notation):\\
{\em Two standard spines of a $3$-manifold $W$ related by a sequence
  of {\rm L}-,~{\rm MP}-,~and false {\rm L}-moves are related by a
  sequence of {\rm L}-~and {\rm MP}-moves only.}\\
By obviously generalizing this proposition to our setting, we obtain a
sequence of L-~and MP-moves only, transforming $R_1$ into $R_2$.
By adding the discs $D^{(j)}$, we obviously obtain the desired
sequence of \La-~and \MPa-moves transforming $(P_1,\beta_1)$ into
$(P_2,\beta_2)$.
{\hfill\hbox{\enspace\fbox{\ref{gener_MP:teo}}}}\vspace{5pt}

\section{Existence of dominating marked spines}\label{dominating:sec}

In this section we generalize, to the setting of marked spines, a
result of Makovetskii~\cite{makov} on the existence of a spine which
dominates, as far as the positive L-moves and positive MP-moves are
concerned, any two given spines of $M$.
Namely, we prove the following.

\begin{teo}\label{gener_makov:teo}
Let $(\calT_1,\beta_1)$ and $(\calT_2,\beta_2)$ be two marked ideal
triangulations of a pair $(M,\alpha)$.
Then there exists a marked ideal triangulation $(\calT,\beta)$ of
$(M,\alpha)$ obtained from both $(\calT_1,\beta_1)$ and
$(\calT_2,\beta_2)$ via a sequence of positive \La-moves and positive
\MPa-moves.
\end{teo}

For the proof, we follow the ideas of~\cite{makov}.

\subsection{Divided spines and moves}

Let us give some definitions useful for the proof.

\paragraph{Dividing strips and divided spines}
Let $(P,\beta)$, with $\beta = \{\beta^{(1)},\ldots ,\beta^{(n)}\}$,
be a spine of a pair $(M,\alpha)$.
Let $\gamma: [0,1] \rightarrow P$ be a simple curve such that:
\begin{itemize}

\item the endpoints belong to edges (maybe, to the same edge) of $P$;

\item $\gamma$ intersects the singularities of $P$ transversely;

\item $\gamma$ contains no vertex of $P$;

\item there exists a strip $S = \gamma \times [0,1] \subset M$
  intersecting $P$ exactly in $\gamma = \gamma \times \{0\}$ and
  $\{\gamma(0),\gamma(1)\} \times [0,1]$.

\end{itemize}
Such a curve $\gamma$ divides some regions of $P$ (those it touches)
into discs: for each region $R$, we will call {\em sub-regions} ({\em
  of $R$}) such discs if $R$ is divided by $\gamma$, or $R$ itself if
it is untouched by $\gamma$.
Let $\overline{\beta} = \{\overline{\beta}^{(1)},\ldots
,\overline{\beta}^{(n)}\}$ be a collection of sub-regions such that
each $\overline{\beta}^{(i)}$ is a sub-region of $\beta^{(i)}$.
The pair $(S,\overline{\beta})$, where $S = \gamma \times [0,1]$, will
be said {\em dividing strip} of $(P,\beta)$, and the triplet
$(P,S,\overline{\beta})$ will be said a {\em divided spine} of
$(M,\alpha)$.

\paragraph{Moves on divided spines}
Let $(P,S,\overline{\beta})$ be a divided spine of $(M,\alpha)$.
We start by defining the obvious generalizations of the positive
\La-~and \MPa-moves and then we define two new moves to take into
account the strip $S$.
As for admissible moves on marked spines, we will say that an
admissible move from $(P,\beta)$ to $(P',\beta')$ gives rise to a {\em
  divided-admissible move} if there is a dividing strip
$(S',\overline{\beta'})$ of $(P',\beta')$ such that
$(P',S',\overline{\beta'})$ is a divided spine of $(M,\alpha)$, and
$(S',\overline{\beta'})$ coincides with $(S,\overline{\beta})$ except
``near'' the portion of $P$ affected by the move.
As it turns out, divided-admissibility depends on $S$.
Moreover, $\overline{\beta'}$ is sometimes not unique.

\subparagraph{MPd-move}
Let us consider a positive \MPa-move $m$ from $(P,\beta)$ to another
spine $(P',\beta')$ of $(M,\alpha)$, such that the strip $S$ is not
involved in the move (namely, $S$ does not intersect the part of $P$
affected by $m$).
Then, we will say that $m$ gives rise to an {\em \MPag-move} from
$(P,S,\overline{\beta})$ to $(P',S',\overline{\beta'})$ whatever
$\overline{\beta}$, where $S'$ coincides with $S$ and
$\overline{\beta'}$ consists of the same sub-regions as
$\overline{\beta}$ (recall that the newborn triangular region does not
belong to $\beta'$).
Note that an \MPag-move always increases (by one) the number of
vertices of $P$.

\subparagraph{Ld-move}
For the \La-moves, the situation is more complicated.
Let us consider a positive \La-move $m$ from $(P,\beta)$ to another
spine $(P',\beta')$ of $(M,\alpha)$, such that the strip $S$ is not
involved in the move (namely, $S$ does not intersect the part of $P$
affected by $m$).
As above, we will say that $m$ gives rise to an {\em \Lag-move} from
$(P,S,\overline{\beta})$ to $(P',S',\overline{\beta'})$ whatever
$\overline{\beta}$, where $S'$ coincides with $S$ and
$\overline{\beta'}$ is uniquely determined by $\overline{\beta}$ and
$\beta'$.
Let us describe $\overline{\beta'}$.
Recall that $m$ divides a region $R$ of $P$ in two regions $R_1$ and
$R_2$, see Fig.~\ref{L_move:fig}-left.
Since the strip $S$ is not involved in the move $m$, then the
\Lag-move divides a sub-region $\overline{R}$ of
$(P,S,\overline{\beta})$ in two sub-regions $\overline{R_1}$ and
$\overline{R_2}$ (where $\overline{R_i}$ is a sub-region of $R_i$, for
$i=1,2$).
Now, we have two cases depending on whether $\overline{R}$ belongs to
$\overline{\beta}$ or not.
If $\overline{R}$ does not belong to $\overline{\beta}$, then
$\overline{\beta'}$ consists of the same sub-regions as
$\overline{\beta}$ ({\em i.e.}~$\overline{R_1}$, $\overline{R_2}$, and
the newborn little region $D$ do not belong to $\beta'$).
If $\overline{R}$ belongs to $\overline{\beta}$, then we define
$\overline{\beta'}$ as $(\overline{\beta} \setminus \{\overline{R}\})
\cup \{\overline{R_1}\}$ or $(\overline{\beta} \setminus
\{\overline{R}\}) \cup \{\overline{R_2}\}$ depending on which region,
between $R_1$ and $R_2$, belongs to $\beta'$.
Note that an \Lag-move always increases (by two) the number of
vertices of $P$.

\subparagraph{Md-move}
We call {\em \Mag-move} any move from a divided spine
$(P,S,\overline{\beta})$ of $(M,\alpha)$ to another divided spine
$(P',S',\overline{\beta'})$ of $(M,\alpha)$, where:
\begin{itemize}

\item $P'$ coincides with $P$;

\item $S'$ is obtained from $S$ as in Fig.~\ref{M_move:fig} (we have
  two cases depending on whether the endpoints of $\gamma$ are
  involved in the move or not); 
  \begin{figure}
    \centerline{\includegraphics{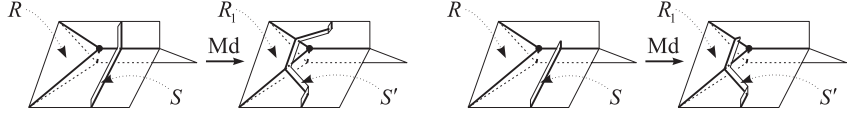}}
    \caption{The \Mag-move (the two cases).}
    \label{M_move:fig}
  \end{figure}
  
\item $\overline{\beta'}$ coincides with $\overline{\beta}$ except
  that the sub-region $R$, if it lies in $\overline{\beta}$, gets
  replaced by the sub-region $R_1$.

\end{itemize}
Note that an \Mag-move increases (by one) the number of intersections
between $\gamma$ and the singularity of $P$, so it can be considered
as being ``positive''.

\subparagraph{Nd-move}
We call {\em \Nag-move} any move from a divided spine
$(P,S,\overline{\beta})$ of $(M,\alpha)$ to another divided spine
$(P',S',\overline{\beta'})$ of $(M,\alpha)$, where:
\begin{itemize}

\item $P'$ coincides with $P$;

\item $S'$ is obtained from $S$ as in Fig.~\ref{N_move:fig};
  \begin{figure}
    \centerline{\includegraphics{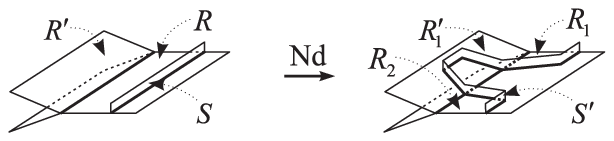}}
    \caption{The \Nag-move.}
    \label{N_move:fig}
  \end{figure}

\item $\overline{\beta'}$ coincides with $\overline{\beta}$ except
  that the sub-regions $R$ and $R'$, if they lie in
  $\overline{\beta}$, get replaced respectively by either the
  sub-region $R_1$ or $R_2$, and by $R'_1$.

\end{itemize}
Note that the choice of which sub-region, between $R_1$ and $R_2$,
belongs to $\overline{\beta'}$ is included in the move.
Finally, note that an \Nag-move increases (by two) the number of
intersections between $\gamma$ and the singularity of $P$, so it can
be considered as being ``positive''.

\subparagraph{} If a spine $(P_2,\beta_2)$ is obtained from a spine
$(P_1,\beta_1)$ via positive \La-moves and positive \MPa-moves, we
will write $(P_1,\beta_1) \nearrow (P_2,\beta_2)$.
If a divided spine $(P_2,S_2,\overline{\beta_2})$ is obtained from a
divided spine $(P_1,S_1,\overline{\beta_1})$ via \Mag-,~\Nag-,~\Lag-,
and \MPag-moves, we will write $(P_1,S_1,\overline{\beta_1}) \nearrow
(P_2,S_2,\overline{\beta_2})$.

\begin{figure}
  \centerline{\includegraphics{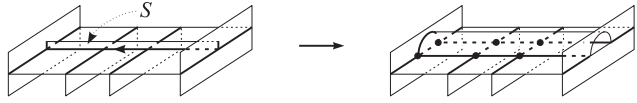}}
  \caption{The swelling.}
  \label{swelling:fig}
\end{figure}

\paragraph{Swelling}
Now we define another move which, taking into account the dividing
strip, transforms a divided spine into a (marked) spine.
Let $(P,S,\overline{\beta})$ be a divided spine of a pair
$(M,\alpha)$, where $S = \gamma \times [0,1]$.
If we apply $m$ positive L-moves to $P$ along the curve $\gamma$
(following the orientation of $\gamma$), we obtain a spine, say $P'$,
of $M$, see Fig.~\ref{swelling:fig}.
Note that $m$ is one less than the number of intersections between
$\gamma$ and the singularities of $P$.
Noting that the collection $\overline{\beta}$ allows us to choose what
regions of $P'$ remain in $\alpha$ after the L-moves (with a little
abuse of notation, we continue calling $\overline{\beta}$ the
collection of such regions), it turns out that the L-moves are
admissible and that the pair $(P',\overline{\beta})$ is a marked spine
of $(M,\alpha)$.
The spine $(P',\overline{\beta})$ will be called {\em swelling of
  $(P,\beta)$ along $(S,\overline{\beta})$} and will be denoted by
$\swell{P,S,\overline{\beta}}$.

\begin{rem}\label{swelling:rem}
For future reference, we underline the fact that\\
$(P,\beta) \nearrow \swell{P,S,\overline{\beta}}$.
\end{rem}

\subsection{Existence of dominating marked spines}

Let us start with two preliminary results.

\begin{figure}
  \centerline{\includegraphics{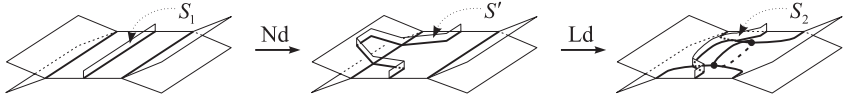}}
  \caption{$(P_1,S_1,\overline{\beta_1}) \nearrow
    (P_2,S_2,\overline{\beta_2})$: the case of a positive \La-move.}
  \label{dividing_curve_lemma:fig}
\end{figure}

\begin{lemma}\label{dividing_curve:lem}
Let $(P_1,S_1,\overline{\beta_1})$ be a divided spine of a pair
$(M,\alpha)$ and let $(P_2,\beta_2)$ be a spine of the pair
$(M,\alpha)$ such that $(P_1,\beta_1) \nearrow (P_2,\beta_2)$.
Then there exists a dividing strip $(S_2,\overline{\beta_2})$ of
$(P_2,\beta_2)$ such that $(P_1,S_1,\overline{\beta_1}) \nearrow
(P_2,S_2,\overline{\beta_2})$.
\end{lemma}
\dimo{dividing_curve:lem}
An easy induction on the number of positive moves transforming
$(P_1,\beta_1)$ into $(P_2,\beta_2)$ allows us to analyze only the
case of only one positive move between $(P_1,\beta_1)$ and
$(P_2,\beta_2)$.
There are two moves to analyze: the positive \La-move and the positive
\MPa-moves.
We concentrate on the first one (the second one being simpler).
If necessary, we first apply \Nag-moves to take the strip $S_1$ away
from the part of $P_1$ affected by the \La-move, see
Fig.~\ref{dividing_curve_lemma:fig}-left.
Let us call $S'$ the strip just obtained.
We impose that the collection $\overline{\beta'}$ consists of the same
sub-regions as $\overline{\beta_1}$, unless a sub-region divided by
one of these \Nag-moves belongs to $\overline{\beta_1}$, in which case
we choose which of the two new sub-regions belongs to
$\overline{\beta'}$ following the choice given by the positive
\La-move.
So $(P_1,S',\overline{\beta'})$ is a divided spine of $(M,\alpha)$.
Now we are able to apply an \Lag-move to $(P_1,S',\overline{\beta'})$
to obtain a divided spine $(P_2,S_2,\overline{\beta_2})$, see
Fig.~\ref{dividing_curve_lemma:fig}-right.
The pair $(S_2,\overline{\beta_2})$ is the dividing strip we are
searching for.
\finedimo{dividing_curve:lem}

\begin{figure}
  \centerline{\includegraphics{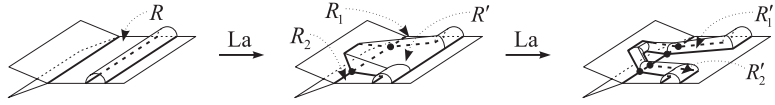}}
  \caption{If $(P_1,S_1,\overline{\beta_1}) \nearrow
    (P_2,S_2,\overline{\beta_2})$ via an \Nag-move, then
    $\swell{P_1,S_1,\overline{\beta_1}} \nearrow
    \swell{P_2,S_2,\overline{\beta_2}}$ via two \La-moves.}
  \label{N_move_swelling_lemma:fig}
\end{figure}

\begin{lemma}\label{swelling:lem}
If $(P_1,S_1,\overline{\beta_1}) \nearrow
(P_2,S_2,\overline{\beta_2})$, then
$\swell{P_1,S_1,\overline{\beta_1}} \nearrow
\swell{P_2,S_2,\overline{\beta_2}}$.
\end{lemma}
\dimo{swelling:lem}
An easy induction on the number of moves transforming
$(P_1,S_1,\overline{\beta_1})$ into $(P_2,S_2,\overline{\beta_2})$
allows us to analyze only the case of only one move between
$(P_1,S_1,\overline{\beta_1})$ and $(P_2,S_2,\overline{\beta_2})$.
There are four possible moves.
If the move is an \Lag-~or an \MPag-move, then obviously
$\swell{P_1,S_1,\overline{\beta_1}} \nearrow
\swell{P_2,S_2,\overline{\beta_2}}$ because $S_1$ is ``far'' from the
move.
If the move is an \Mag-move, we have three cases:
\begin{enumerate}

\item If $\gamma(0)$ is involved in the move (see
  Fig.~\ref{M_move:fig}-right), then
  $\swell{P_2,S_2,\overline{\beta_2}}$ is obtained from
  $\swell{P_1,S_1,\overline{\beta_1}}$ via a positive \La-move, as
  shown in Fig.~\ref{M_move_swelling_lemma:fig}.
  \begin{figure}
    \centerline{\includegraphics{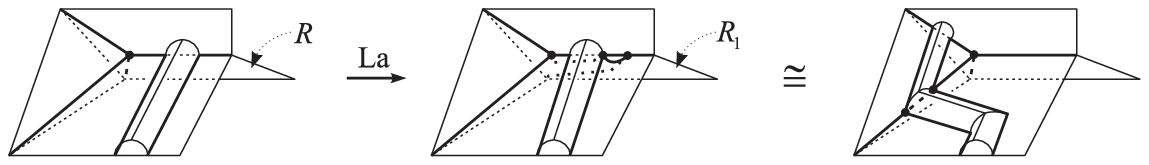}}
    \caption{If $(P_1,S_1,\overline{\beta_1}) \nearrow
      (P_2,S_2,\overline{\beta_2})$ via an \Mag-move, then
      $\swell{P_1,S_1,\overline{\beta_1}} \nearrow
      \swell{P_2,S_2,\overline{\beta_2}}$ via an \La-move (case~1).}
    \label{M_move_swelling_lemma:fig}
  \end{figure}
  Note that, if the region $R$ belongs to $\overline{\beta_1}$, we
  choose to leave in $\overline{\beta_2}$ the region $R_1$; so the
  spine obtained is exactly the swelling of $(P_2,\beta_2)$ along
  $(S_2,\overline{\beta_2})$.

\item If neither $\gamma(0)$ nor $\gamma(1)$ is involved in the move
  (see Fig.~\ref{M_move:fig}-left), then
  $\swell{P_2,S_2,\overline{\beta_2}}$ is obtained from
  $\swell{P_1,S_1,\overline{\beta_1}}$ via two positive \MPa-moves.

\item If $\gamma(1)$ is involved in the move (see again
  Fig.~\ref{M_move:fig}-right), then
  $\swell{P_2,S_2,\overline{\beta_2}}$ is obtained from
  $\swell{P_1,S_1,\overline{\beta_1}}$ via two positive \MPa-moves.

\end{enumerate}
If the move is an \Nag-move (see Fig.~\ref{N_move:fig}), then
$\swell{P_2,S_2,\overline{\beta_2}}$ is obtained from
$\swell{P_1,S_1,\overline{\beta_1}}$ via two positive \La-moves, as
shown in Fig.~\ref{N_move_swelling_lemma:fig}.
For the first \La-move, if the region $R$ belongs to
$\overline{\beta_1}$, we have to choose a region, between $R_1$ and
$R_2$, to leave in $\overline{\beta_2}$: we choose the region
depending on which sub-region, between $R_1$ and $R_2$, belongs to
$\overline{\beta_2}$ after the \Nag-move.
For the second \La-move, if the region $R'$ belongs to
$\overline{\beta_1}$, we choose to leave in $\overline{\beta_2}$ the
``nearest'' region (between $R'_1$ and $R'_2$) to $\gamma(0)$.
The spine obtained is exactly the swelling of $(P_2,\beta_2)$ along
$(S_2,\overline{\beta_2})$.
This concludes the proof.
\finedimo{swelling:lem}

Now we are able to prove Theorem~\ref{gener_makov:teo}.

\dimo{gener_makov:teo}
Let $(P_i,\beta_i)$ the dual spine of $(\calT_i,\beta_i)$, for
$i=1,2$.
By applying Theorem~\ref{gener_MP:teo} and by noting that each
\Va-move is actually an \La-move, we obtain a sequence $s$ of \La-~and
\MPa-moves transforming $(\calT_1,\beta_1)$ into $(\calT_2,\beta_2)$.
The sequence $s$ can be divided in (sub-)sequences $s_i$, with
$i=1,\ldots ,2l$, where the sequences $s_{2k+1}$ are composed by
positive moves while the sequences $s_{2k}$ are composed by negative
moves, and only $s_1$ and $s_{2l}$ could be empty.
Let us call $|s_i|$ the number of moves of the sequence $s_i$.
An easy induction on $S=\sum_{k=1}^{l-1} |s_{2k+1}|$ allows us to
prove only the following statement.

{\em If $(P_2,\beta_2)$ is obtained from $(P_1,\beta_1)$ via a
  sequence $s$ such that $l=2$, $|s_1|=0$, $|s_2|>0$, $|s_3|=1$ and
  $|s_4|=0$, then there exists another sequence $s'$, transforming
  $(P_1,\beta_1)$ into $(P_2,\beta_2)$, such that $l=1$}.

The proof of this statement concludes the proof of the theorem.
We have to prove that there exists a spine $(P,\beta)$ such that
$(P_1,\beta_1) \nearrow (P,\beta) \nwarrow (P_2,\beta_2)$.
If we call $(Q,\beta')$ the spine before the positive move $m$ of the
sequence $s_3$, we have that $(P_1,\beta_1) \nwarrow (Q,\beta')
\nearrow (P_2,\beta_2)$.
Let us start by choosing a dividing strip $(S',\overline{\beta'})$ for
$(Q,\beta')$ (we have two cases depending on $m$).
\begin{itemize}

\item If $m$ is a positive \La-move, we choose as $\gamma'$ the curve
  determining $m$.
  Note that there are two different strips $S' = \gamma' \times [0,1]$
  (up to isotopy): we choose one of them (the choice is immaterial).
  If the region of $Q$ divided by $\gamma'$ is one of the
  $(\beta')^{(i)}$'s, we choose the $(\overline{\beta'})^{(i)}$
  following the choice given by $m$.

\item If $m$ is a positive \MPa-move, we choose as $\gamma'$ a curve
  parallel to the edge $e$ of $Q$ disappearing during $m$.
  As above there are two different strips: we choose one of them.
  If the region of $Q$ divided by $\gamma'$ is one of the
  $(\beta')^{(i)}$'s, we choose as $(\overline{\beta'})^{(i)}$ the
  sub-region which is not adjacent (locally) to $e$.

\end{itemize}

By Lemma~\ref{dividing_curve:lem}, there exists a dividing strip
$(S_1,\overline{\beta_1})$ for $(P_1,\beta_1)$ such that
$(P_1,S_1,\overline{\beta_1}) \nwarrow (Q,S',\overline{\beta'})$; so,
by Lemma~\ref{swelling:lem}, $\swell{P_1,S_1,\overline{\beta_1}}
\nwarrow \swell{Q,S',\overline{\beta'}}$.
By Remark~\ref{swelling:rem}, we have that $(P_1,\beta_1) \nearrow
\swell{P_1,S_1,\overline{\beta_1}}$.
Finally, we have two cases depending on $m$.
\begin{itemize}

\item If $m$ is a positive \La-move, then
  $\swell{Q,S',\overline{\beta'}} = (P_2,\beta_2)$; so we have that
  $$
  (P_1,\beta_1) \nearrow \swell{P_1,S_1,\overline{\beta_1}} \nwarrow
  \swell{Q,S',\overline{\beta'}} = (P_2,\beta_2).
  $$

\item If $m$ is a positive \MPa-move, then
  $\swell{Q,S',\overline{\beta'}}$ can be obtained from
  $(P_2,\beta_2)$ via a positive \MPa-move, see
  Fig.~\ref{MP_swelling:fig}; so we have that
  $$
  (P_1,\beta_1) \nearrow \swell{P_1,S_1,\overline{\beta_1}} \nwarrow
  \swell{Q,S',\overline{\beta'}} \nwarrow (P_2,\beta_2).
  $$
  \begin{figure}
    \centerline{\includegraphics{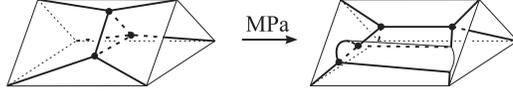}}
    \caption{If $m$ is a positive \MPa-move, then $(P_2,\beta_2)
      \nearrow \swell{Q,S',\overline{\beta'}}$.}
    \label{MP_swelling:fig}
  \end{figure}

\end{itemize}
\finedimo{gener_makov:teo}

\section{Applications}

In this section we describe two applications of the previous results.
The first one is due to Basehilac and Benedetti~\cite{BB1,BB2,BB3}.
The second one is a natural question arisen in a work of Frigerio and
Petronio~\cite{Frig-Petr}.

\subsection{Links in 3-manifolds}\label{hamiltonian:subsec}

Let $M$ be a closed 3-manifold and $L$ a link in $M$.
A pair $(\calT,\calL)$ is said to be a {\em distinguished
  triangulation} of the pair $(M,L)$ if $\calT$ is a loose
triangulation of $M$, the link $L$ is triangulated by $\calL$ and
$\calL$ is a Hamiltonian sub-complex of $\calT$ ({\em i.e.}~each
vertex of $\calT$ is an endpoint of exactly two germs of edges of
$\calL$).
As we have done for marked ideal triangulations, we can define ({\em
  positive} and {\em negative}) {\em admissible} MP-~and L-moves
between distinguished triangulations.
We need another move allowing us to change the number of vertices of
$\calT$.
We will say that the distinguished triangulation $(\calT',\calL')$ is
obtained from the distinguished triangulation $(\calT,\calL)$ via a
{\em positive admissible} B-move if
\begin{itemize}
\item
  $\calT'$ is obtained from $\calT$ via a positive B-move,
\item
  one edge $e$ of the tetrahedron $T$ involved in the move belongs to $\calL$,
\item
  $\calL'$ coincides with $\calL$ except for the edge $e$ which is
  substituted with the other two edges of the only triangle of
  $\calT'$ created by the B-move and containing $e$.
  \begin{figure}
    \centerline{\includegraphics{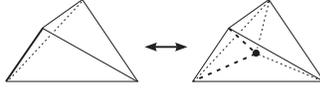}}
    \caption{An admissible B-move on a distinguished triangulation.
      (The link is drawn thick.)}
    \label{distinguished_B_move:fig}
  \end{figure}
\end{itemize}
See Fig.~\ref{distinguished_B_move:fig}.
Obviously, a {\em negative admissible} B-move between distinguished
triangulations is defined as the inverse of a positive admissible
B-move.

Now we are able to prove the calculus for distinguished
triangulations.
\begin{cor}\label{distinguished_calculus:cor}
Two distinguished triangulations of a pair $(M,L)$ can be obtained
from each other via a sequence of admissible {\rm B}-~and
{\rm MP}-moves.
\end{cor}
\dimo{distinguished_calculus:cor}
Let $(\calT_1,\calL_1)$ and $(\calT_2,\calL_2)$ be two distinguished
triangulations of $(M,L)$.
Obviously, up to applying suitable admissible B-moves, we can suppose
that $(\calT_1,\calL_1)$ and $(\calT_2,\calL_2)$ have the same number
of vertices on each component of $L$.
Moreover, up to isotopy, we can suppose that the links $\calL_i$
coincide with $L$, and that the vertices of $\calT_1$ and the vertices
of $\calT_2$ coincide with each other.

Now, we remove a little star of each vertex of $\calT_i$: let us call
$\overline{M}$ the manifold just obtained.
Obviously, after removing the balls, the link $L$ becomes a collection
of arcs, say $\overline{L}$, and, for $i=1,2$, the pair
$(\calT_i,\calL_i)$ is a marked loose triangulation corresponding to a
marked ideal triangulation of $(\overline{M},\overline{L})$.
So, by applying Corollary~\ref{gener_MP_senza_V:cor}, we obtain that
$(\calT_2,\calL_2)$ can be obtained from $(\calT_1,\calL_1)$ via
admissible MP-moves.
This concludes the proof.
\finedimo{distinguished_calculus:cor}

Using the same technique (and Theorem~\ref{gener_makov:teo}), the
following result on dominating distinguished triangulations can be
proved.
\begin{cor}
Let $(\calT_1,\calL_1)$ and $(\calT_2,\calL_2)$ be two distinguished
triangulations of a pair $(M,L)$.
Then there exists a distinguished triangulation $(\calT,\calL)$ of
$(M,L)$ obtained from both $(\calT_1,\calL_1)$ and $(\calT_2,\calL_2)$
via a sequence of admissible positive {\rm B}-,~{\rm L}-,~and {\rm
  MP}-moves.
\end{cor}

\subsection{Partially truncated triangulations}\label{part_trunc_tria:subsec}

In this subsection we briefly describe a generalization of ideal
triangulations which is useful to study complete finite-volume
orientable hyperbolic 3-manifolds with geodesic
boundary~\cite{Frig-Petr}.
(For the sake of shortness, in the rest of the subsection we will just
say {\em hyperbolic}.)
For a complete description see~\cite{Frig-Petr}.

Let $N$ be such a hyperbolic manifold.
It is a fact that $N$ consists of a compact portion together with some
cusps based either on tori or on annuli.
This fact implies that $N$ has a natural compactification
$\overline{N}$ obtained from $N$ by adding some tori and some annuli.
Let us call $C$ and $A$ the collection of such tori and such annuli,
respectively.
It is a fact that $N$ can be obtained in a non-ambiguous way from the
pair $(\overline{N},A)$ by removing from $\overline{N}$ both $A$ and
{\em all} the toric components of $\partial\overline{N}$.
Moreover, there is no sphere in $\partial\overline{N}$ and there is no
annulus in $A$ which is contained in a torus of $C$.

\begin{figure}
  \centerline{\includegraphics{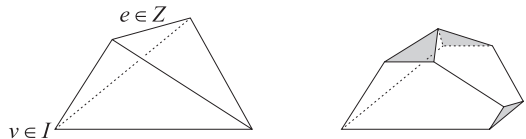}}
  \caption{A partially truncated tetrahedron with one ideal vertex and
    one length-0 edge (on the left) and its topological realization
    (on the right).}
  \label{part_trunc_tetra:fig}
\end{figure}

Let us describe now a generalization of ideal triangulations, which
takes into account the annuli.
Let us start by defining the pieces substituting ideal tetrahedra.
A {\em partially truncated tetrahedron} is a triple $(T,I,Z)$ where
$T$ is a tetrahedron, $I$ is a set of vertices of $T$ (called {\em
  ideal vertices}), and $Z$ is a set of edges of $T$ (called {\em
  length-$0$ edges}) such that neither of the two endpoints of an edge
in $Z$ belongs to $I$.
Now we define the {\em topological realization} of a partially
truncated tetrahedron $(T,I,Z)$ as the space $T^*$ obtained by
removing from the tetrahedron $T$ the ideal vertices, the length-0
edges, and small open stars of the non-ideal vertices.
We call {\em lateral hexagon} and {\em truncation triangle} the
intersection of $T^*$ respectively with a face of $T$ and with the
link in $T$ of a non-ideal vertex.
Note that, if $(T,I,Z)$ has a length-0 edge, some vertices of a
truncation triangle of $T^*$ may be missing and, if $(T,I,Z)$ has
ideal vertices or length-0 edges, a lateral hexagon of $T^*$ may not
be a hexagon, because some of its edges may be missing.
See Fig.~\ref{part_trunc_tetra:fig}.

Let us consider now a manifold $N$ which is a candidate to be
hyperbolic.
Namely, let $\overline{N}$ be a compact orientable manifold, having no
sphere in the boundary, and let $A \subset \partial\overline{N}$ be a
family of disjoint annuli not lying on the toric components of
$\partial\overline{N}$; let $N$ be obtained from $\overline{N}$ by
removing $A$ and the toric components.
Finally, we define a {\em partially truncated triangulation} of $N$ as
a realization of $N$ as the gluing of some $T^*$'s along a pairing of
the lateral hexagons induced by a simplicial pairing of the faces of
the $T$'s.
Note that the truncation triangles of the $T^*$'s give a triangulation
of $\partial N$ with some genuine and some ideal vertices, the links
of the ideal vertices of the $T$'s give a triangulation of the toric
components of $\partial\overline{N}$, and the links of the length-0
edges of the $T$'s give a decomposition into rectangles of the annuli
in $A$.

Let us now translate the theory of partially truncated triangulations
into the language of marked ideal triangulations.
Let us consider $\overline{N}$ as above and let us collapse every
annulus $[-1,1]\times S^1 \in A$ to an arc $[-1,1]\times\{*\}$.
It turns out that the space just obtained, say $N'$, is a compact
3-manifold and each $[-1,1]\times\{*\}$ is an arc properly embedded in
$N'$.
Let us call $\alpha_N$ the family of the arcs $[-1,1]\times\{*\}$ in
$N'$.
It is a fact that partially truncated triangulations of $N$
bijectively correspond to marked ideal triangulations of the pair
$(N',\alpha_N)$; under this correspondence, the length-0 edges and the
ideal vertices correspond respectively to the edges in $\alpha_N$ and
to the vertices on the tori of $\partial N'$ on which there are no
ends of arcs in $\alpha_N$.

Obviously, the admissible MP-~and V-moves between marked ideal
triangulations of $(N,\alpha_N)$ translate into moves between
partially truncated triangulations of $N$.
Let us call {\em admissible {\rm MP}-~{\rm and V}-moves} also such
moves between partially truncated triangulations.
Now, Theorem~\ref{gener_MP:teo} and
Corollary~\ref{gener_MP_senza_V:cor} imply the following.
\begin{cor}
Two partially truncated triangulations of $N$ can be obtained from
each other via a sequence of admissible {\rm V}-~and {\rm MP}-moves.
If moreover the two partially truncated triangulations have at least
two tetrahedra, then they can be obtained from each other via a
sequence of admissible {\rm MP}-moves only.
\end{cor}

\vspace{0.5cm}

\noindent amendola@mail.dm.unipi.it,\\ 
Dipartimento di Matematica,\\ 
Via F. Buonarroti 2,\\ 
I-56127 PISA


\begin{thebibliography}{99}



\bibitem{Alexander}
\textsc{J.~W.~Alexander},
\textit{The combinatorial theory of complexes},
Ann. Math. {\bf 31} (1930), 294-322.


\bibitem{BB1}
\textsc{S.~Basehilac -- R.~Benedetti},
\textit{QHI Theory, {\rm I:} $3$-manifolds scissors congruence
  classes and quantum hyperbolic invariants},
\texttt{math.GT/0201240}.


\bibitem{BB2}
\textsc{S.~Basehilac -- R.~Benedetti},
\textit{QHI, $3$-manifolds scissors classes and the volume
  conjecture},
Geometry and Topology Monographs, Vol.~4: Invariants of Knots and
3-manifolds (Kyoto 2001), 13-28.


\bibitem{BB3}
\textsc{S.~Basehilac -- R.~Benedetti},
\textit{{\rm QHI} Theory, {\rm II}: dilogarithmic and quantum
  hyperbolic invariants of $3$-manifolds with
  $PSL(2,\mathbb{C})$-characters {\rm [Abridged Version]}},
\texttt{math.GT/0211061}.


\bibitem{Casler}
\textsc{B.~G.~Casler},
\textit{An imbedding theorem for connected $3$-manifolds with
  boundary},
Proc. Amer. Math. Soc. {\bf 16} (1965), 559-566.


\bibitem{Frig-Petr}
\textsc{R.~Frigerio -- C.~Petronio},
\textit{Construction and recognition of hyperbolic $3$-manifolds with
  geodesic boundary},
\texttt{math.GT/0109012}.


\bibitem{makov}
\textsc{A.~Yu.~Makovetskii},
\textit{On transformations of special spines and special polyhedra},
Math. Notes {\bf 65} (1999), 295-301.


\bibitem{Matveev:calculus}
\textsc{S.~V.~Matveev},
\textit{Transformations of special spines and the Zeeman conjecture},
Math. USSR-Izv. {\bf 31} (1988), 423-434.


\bibitem{Matveev:new:book}
\textsc{S.~V.~Matveev},
``Algorithmic methods in 3-manifolds topology'',
Springer-Verlag, to appear.


\bibitem{Matveev-Fomenko}
\textsc{S.~V.~Matveev -- A.~T.~Fomenko},
\textit{Constant energy surfaces of Hamiltonian systems, enumeration
  of three-dimensional manifolds in increasing order of complexity,
  and computation of volumes of closed hyperbolic manifolds},
Russ. Math. Surv. {\bf 43} (1988), 3-25.


\bibitem{Pachner}
\textsc{U.~Pachner}.
\textit{P.L. homeomorphic manifolds are equivalent by elementary
  shellings},
Europ. J. Combinatorics {\bf 12} (1991), 129-145.


\bibitem{Petronio:tesi}
\textsc{C.~Petronio}, 
``Standard spines and 3-manifolds'',
Scuola Normale Superiore, Pisa, 1995.


\bibitem{Piergallini}
\textsc{R.~Piergallini},
\textit{Standard moves for standard polyhedra and spines},
Rendiconti Circ. Mat. Palermo {\bf 37}, suppl.~18 (1988), 391-414.


\bibitem{Turaev-Viro}
\textsc{V.~G.~Turaev -- O.~Ya.~Viro},
\textit{State sum invariants of $3$-manifolds and quantum
  $6j$-symbols},
Topology {\bf 31} (1992), 865-902.



\end{thebibliography}
\end{document}